\documentclass[final,onefignum,onetabnum]{siamart250211}

\usepackage{amsmath}
\usepackage{amsfonts}
\usepackage{dsfont}
\usepackage{graphicx}
\usepackage{epstopdf}
\usepackage{lineno}

\ifpdf
\DeclareGraphicsExtensions{.eps,.pdf,.png,.jpg}
\else
\DeclareGraphicsExtensions{.eps}
\fi

\newtheorem{prop}{Property}

\usepackage{xcolor}
\definecolor{darkgreen}{rgb}{0.1,0.6,0.1}

\allowdisplaybreaks

\newsiamremark{remark}{Remark}
\newsiamremark{hypothesis}{Hypothesis}
\crefname{hypothesis}{Hypothesis}{Hypotheses}
\newsiamthm{claim}{Claim}

\headers{Structure preserving FV schemes on Voronoi grids}{W. Boscheri and F. Dhaouadi}    %

\title{Structure preserving finite volume schemes on Voronoi grids: curl involution, asymptotic limit and thermodynamics}

\author{Walter Boscheri \thanks{Laboratoire de Mathématiques UMR 5127 CNRS, Universit{\'e} Savoie Mont Blanc, 73376 Le Bourget du Lac, France (\email{walter.boscheri@univ-smb.fr})} 
	\and 	
	Firas Dhaouadi\thanks{Department of Civil, Environmental and Mechanical Engineering, University of Trento, Via Mesiano 77, 38123 Trento, Italy (\email{firas.dhaouadi@unitn.it})} 
}

\newcommand{\prn}[1]{\left( #1 \right)} 
\newcommand{\pd}[2]{\frac{\partial #1}{\partial #2}}
\newcommand{\Nc}{N_c}
\newcommand{\x}{\mathbf{x}}
\newcommand{\n}{\mathbf{n}}
\renewcommand{\t}{\mathbf{t}}
\newcommand{\norm}[1]{\left\vert \left\vert #1 \right\vert \right \vert}
\newcommand{\abs}[1]{\left\vert  #1 \right \vert}
\newcommand{\ppsi}{\boldsymbol{\psi}}

\renewcommand{\j}{\mathbf{j}}
\newcommand{\temp}{\theta}        
\newcommand{\rhos}{s}             
\newcommand{\vv}{\mathbf{u}}      
\newcommand{\p}{\mathbf{p}}       
\newcommand{\w}{\mathbf{w}}       
\newcommand{\dt}{\Delta t}        
\newcommand{\ar}{\alpha(\rho)}    
\newcommand{\f}{\mathbf{f}}       
\newcommand{\g}{\mathbf{g}}       
\newcommand{\h}{\mathbf{h}}       
\newcommand{\as}{\delta^{\rhos}_f}  
\newcommand{\A}{\mathbf{A}}       
\newcommand{\J}{\mathbf{J}}       
\newcommand{\D}{\mathbf{D}}       
\newcommand{\B}{\mathbf{B}}       
\newcommand{\F}{\mathbf{F}}       
\renewcommand{\L}{\mathbf{L}}     

\begin{document}

\maketitle

\begin{abstract}
We propose a new curl-free and thermodynamically compatible finite volume scheme on Voronoi grids to solve compressible heat conducting flows written in first-order hyperbolic form. The approach is based on the definition of compatible discrete curl-grad operators, exploiting the triangular nature of the dual mesh. We design a cell solver reminiscent of the nodal solvers used in Lagrangian schemes to discretize the evolution equation for the thermal impulse vector, and we demonstrate that the resulting numerical scheme ensures energy conservation, local non-negative entropy production, as well as asymptotic consistency with the classical Fourier law in the stiff relaxation limit. A novel technique is proposed to transfer residuals from the dual to the primal mesh as subfluxes, which eventually yields the construction of entropy compatible semi-discrete methods. The scheme and its properties are validated on a set of numerical test cases.  
\end{abstract}



\section{Introduction} 

The numerical approximation of nonlinear hyperbolic systems poses major challenges, as solutions may exhibit discontinuities even for smooth initial data. Selecting the physically relevant weak solution requires that numerical methods abide by the principles of thermodynamics and in particular energy conservation (first principle) and non-negative entropy production (second principle). A classical example is the Euler equations of gas dynamics, where mass, momentum and energy must be conserved while ensuring entropy increase across shocks \cite{godunov1959finite,ToroBook}. However, the consistency with physics does not stop at the enforcement of the principles of thermodynamics. Indeed, many models of continuum mechanics also include geometric or structural constraints, such as divergence-free and curl-free involutions, that must be preserved over time. One can recall for instance the Maxwell equations of electrodynamics \cite{maxwell1865viii}, magnetohydrodynamics \cite{alfven1942existence} or the Einstein field equations \cite{einstein1916foundation}. Failure to preserve these constraints at the discrete level may lead to a degradation of the solution quality, finite-time blowup \cite{chiocchetti2023exactly, dhaouadi2023structure} or violation of fundamental physical properties \cite{brackbill1980effect,brackbill1985fluid,balsara1999staggered}. From a modeling perspective, relaxation systems \cite{jin1995relaxation} have proven to be powerful tools for constructing first-order hyperbolic models \cite{PeshRom2014}. In the stiff relaxation limit, the original governing equations are retrieved, which is an additional property that must be preserved also at the discrete level. Asymptotic preserving methods \cite{filbet2010} ensure this consistency with the underlying physical modeling.

One particular model in which the aforementioned complications rise is the system of inviscid heat conducting compressible fluids introduced in \cite{dhaouadi2024eulerian}. This system extends the Euler equations for compressible flow to include heat transfer, using an unconditionally hyperbolic model consistent with the first and second laws of thermodynamics, instead of the traditional parabolic heat flux. The description of heat propagation relies on an additional vector field denoted as $\j$, which is asymptotically related to the Fourier heat flux, when the relaxation time goes to zero. Hyperbolic modeling of heat conduction began with Cattaneo \cite{cattaneo1948}, whose model was later made Galilean-invariant \cite{christov2009frame}. However, its coupling with compressible flows introduces non-hyperbolicity in multiple dimensions and thermodynamic inconsistencies \cite{angeles2022non,rubin1992hyperbolic}. The approach proposed in \cite{peshkov2018continuum} provides a sounder alternative that is thermodynamically compatible with hyperbolicity established in the vicinity of equilibrium states.

In order to solve the system of heat conducting compressible fluids \cite{dhaouadi2024eulerian}, we propose in this paper a novel finite volume scheme on Voronoi grids. The scheme relies on a dual grid structure such that the curl-free vector $\j$ is stored and evolved in the vertices of the Voronoi grid using conservative compatible gradient operators, while the rest of the variables are evolved in the Voronoi cell centers using a classical finite volume method. The scheme is designed such that the following properties hold at the discrete level.
\begin{enumerate}
	\item \textbf{Curl-involution preservation:} In the absence of relaxation source term, the scheme preserves exactly the curl involution constraint on the vector $\j$, inside the discrete domain as well as on its boundaries.
	\item  \textbf{Asymptotic preservation:} The scheme preserves the asymptotic limit of the model in which the parabolic Fourier law is recovered.
	\item \textbf{Thermodynamic compatibility:} The discrete total energy is conserved, and the discrete entropy balance is conserved for smooth flows. 
\end{enumerate}
The first property is achieved by modifying the curl-grad operators proposed in \cite{barsukow2025node} to deal with Voronoi grid and its associated dual triangular mesh. The second property is ensured by using an implicit-explicit scheme \cite{pareschirusso2005} that retrieves the Fourier law in the first-order asymptotic expansion of $\j$. The third property is less straightforward since we are dealing with vertex-staggered meshes: entropy production is guaranteed by a novel cell solver for the numerical flux of the $\j$ equation inspired by \cite{maire2011staggered}, while entropy balance is respected for smooth flows with an innovative technique that transfers the residuals from the dual to the primal mesh, exploiting a subflux reconstruction along the lines of \cite{Abgrall_subflux2018}. Thermodynamic compatibility is then restored at the semi-discrete level by means of a scalar correction factor \cite{Abgrall2018}. To the best of our knowledge, these discrete properties have never been achieved simultaneously on polygonal grids, and our framework can be easily extended to other type of hyperbolic models.

This paper is organized as follows. In Sec.~\ref{sec:Model}, we recall the governing system of equations and its main properties, namely its hyperbolicity, eigenstructure and asymptotic behavior. In Sec.~\ref{sec:NumScheme} we thoroughly describe the numerical method and demonstrate its properties. In particular, we define the dual mesh configuration and related differential operators from which the discrete curl-grad identity follows as a direct consequence. A curl-free discretization is then constructed and analyzed for the vector $\j$ and it is shown to guarantee a positive discrete entropy production and proven to be asymptotically consistent with the Fourier limit. The novel technique to transfer residuals of $\j$ from the dual to the primal mesh is then detailed, proving entropy preservation on the main grid for smooth flows. 
In Sec.~\ref{sec:Results}, we validate the proposed numerical scheme on a set of designed test cases and benchmarks. Finally, we conclude the paper and give an outlook to future research lines in Sec.~\ref{sec:Conclusions}.
\section{Hyperbolic heat transfer in compressible flows}
\label{sec:Model}

\subsection{Governing equations}

The system of equations proposed in \cite{dhaouadi2024eulerian} and describing the motion of heat conducting compressible flows is given by 
\begin{subequations}
	\begin{alignat}{2}
		&\pd{\rho}{t} &&+ \nabla \cdot (\rho \mathbf{u}) = 0, \label{eqn.rho} \\
		&\pd{\rho \mathbf{u}}{t} &&+ \nabla \cdot \left( \rho \mathbf{u} \otimes \mathbf{u} + \left( p + \tfrac{1}{2} (\rho \alpha'(\rho) - \ar) \norm{\mathbf{j}}^2 \right) \mathbf{I} + \ar \, \mathbf{j} \otimes \mathbf{j} \right) = \mathbf{0}, \label{eqn.momentum} \\
		&\frac{\partial \j}{\partial t} &&+ \nabla\prn{\j\cdot\mathbf{u} + \theta(\rho,\eta)}  + \prn{ \frac{\partial \j}{\partial \x} - \prn{\frac{\partial \j}{\partial \x}}^T }\mathbf{u} =  -\frac{\mathbf{j}}{\tau} \label{eqn.J},  \\
		&\pd{E}{t} &&+ \nabla \cdot \left( \left( E + p + \tfrac{1}{2} (\rho \alpha'(\rho) - \ar) \norm{\mathbf{j}}^2 \right) \mathbf{u} + \ar (\mathbf{j} \cdot \mathbf{u}) \mathbf{j} + \mathbf{q} \right) = 0. \label{eqn.en}
	\end{alignat}
	\label{eqn.heat_system}
\end{subequations}
Here, the spatial differential operators are taken with respect to $\mathbf{x}\in \mathbb{R}^d$, where $d=2$ is the number of space dimensions and $t\in \mathbb{R}_+$ is the time. $\rho$~is the fluid density and $\mathbf{u}\in \mathbb{R}^d$ is its velocity field. $\j\in\mathbb{R}^d$ is the so-called thermal impulse and $E$ is the total energy density, expressed in terms of the system variables and the entropy $\eta$ as
\begin{equation}
	E = \frac{1}{2} \rho \norm{\vv}^2 + \rho \varepsilon(\rho,\eta) + \frac{1}{2}\ar \norm{\j}^2,
	\label{eqn.rhoE}
\end{equation}  
where $\ar$ is a positive function. The fluid pressure $p$ and its temperature $\temp$ are related to the system internal energy $\varepsilon(\rho,\eta)$ via the Gibbs identity 
\begin{equation}
	\text{d}\varepsilon + p\, \text{d}\prn{1/\rho}  - \temp\, \text{d}\eta = 0,
	\label{eqn.Gibbs}
\end{equation}
so that 
\begin{equation*}
	p(\rho,\eta) = \rho^2 \pd{\varepsilon(\rho,\eta)}{\rho}, \quad \temp(\rho,\eta) = \pd{\varepsilon(\rho,\eta)}{\eta}.
\end{equation*}
The flux vector $\mathbf{q} = \ar \temp\, \j$, appearing in the energy flux, is due to interstitial working, and in this particular case it can be regarded as the heat flux density vector.
The entropy equation writes 
\begin{equation}
	\pd{\rhos}{t}  + \nabla\cdot \prn{s \mathbf{u} + \ar\,\j } = \frac{\ar}{ \theta\,\tau}\norm{\j}^2 \geq 0, \label{eqn.eta}
\end{equation}
where $s = \rho \eta$. It can be derived by expressing the energy $E$ as a function of the conserved variables $\tilde{\w} = (\rho,\rho\vv,\j,s)$, differentiating with respect to time and deducing \eqref{eqn.eta} from the equality 
\begin{equation}
	\pd{E}{t} = \tilde{\p} \cdot \pd{\tilde{\w}}{t} = \frac{\partial E}{\partial \rho} \pd{\rho}{t} + \frac{\partial E}{\partial \rho \mathbf{u}} \pd{\rho \vv}{t} + \frac{\partial E}{\partial \j} \pd{\j}{t} + \frac{\partial E}{\partial \rhos} \pd{\rhos}{t}  ,
	\label{eqn.pdq}
\end{equation}
where $\tilde{\p}=(\partial_{\rho}E, \, \partial_{\rho \mathbf{u}}E, \, \partial_{\rhos}E, \, \partial_{\j}E)^\top$ is the vector of the conjugate variables. 
\newline In what follows, we will consider an ideal gas equation of state:
\begin{equation}
	\varepsilon(\rho,\eta) = \frac{\rho^{\gamma-1}}{\gamma-1}e^{\eta/c_v}, \quad p(\rho,\eta ) = \rho^\gamma e^{\eta/c_v}, \quad \temp(\rho,\eta) = \frac{\rho^{\gamma-1}}{c_v(\gamma-1)}e^{\eta/c_v},
\end{equation} 
where $c_v$ is the heat capacity at constant volume, and $\gamma$ is the heat capacity ratio.
\subsection{Structural properties of the model}
In the absence of relaxation term, the thermal impulse is subject to a curl-free involution: 
\begin{equation}
\text{If at } t=0, \quad \nabla \times \j = 0 \quad \forall \mathbf{x} \in \Omega \quad \text{then} \quad \nabla \times \j = 0 \quad \forall \mathbf{x} \in \Omega, \ \forall t>0.
\label{eqn.inv}
\end{equation}
The set of equations \eqref{eqn.heat_system} allows to describe heat transfer processes in a compressible flow via a system of first-order Partial Differential Equations (PDE) with stiff source terms. Indeed, a formal expansion in the parameter $\tau$ shows that for the value of
\begin{equation}
	\tau = \frac{K}{\ar \, \temp},
	\label{eqn.tau}
\end{equation}
one recovers the classical Fourier heat flux with heat conductivity $K$ \cite{dhaouadi2024eulerian}:
\begin{equation}
	\mathbf{q} = -\tau \ar 	\, \temp \, \nabla \temp = - K \nabla \temp.
	\label{eqn.heatflux}
\end{equation}
For $d\geq2$, the system of equations \eqref{eqn.heat_system} is weakly hyperbolic. The characteristic speeds $\boldsymbol{\lambda}$ in the $x$-direction can be computed explicitly for the particular case $\ar=~\kappa^2/\rho$ and are given by \cite{dhaouadi2024eulerian}  
\begin{equation}
	\lambda_{1-4} = u_1 \pm \sqrt{Z_1 \pm Z_2}, \quad \lambda_{5-2d+2} = u_1,
\end{equation}	
where
\begin{eqnarray*}
	Z_1 &=& \frac{1}{2}\prn{\gamma \rho^{\gamma-1}e^{\eta/c_v} +\frac{\kappa^2 \rho^{\gamma-3}e^{\eta/c_v}}{c_v^2(\gamma-1)} + 2\frac{\kappa^2}{\rho^2} \prn{J_2^2+J_3^2}}, \label{eqn.lambda}\\
	Z_2 &=& \sqrt{\prn{\frac{\kappa}{c_v} \rho^{\gamma-2}e^{\eta/c_v}}^2 + \frac{1}{4} \prn{\gamma \rho^{\gamma-1}e^{\eta/c_v} -\frac{\kappa^2 \rho^{\gamma-3}e^{\eta/c_v}}{c_v^2(\gamma-1)}}^2}.
\end{eqnarray*}
The weak hyperbolicity is due to the multiple eigenvalue $u_1$ having a degenerate eigenspace, with $d-1$ missing eigenvectors. Nevertheless, strong hyperbolicity can be restored either by supplying the so-called Godunov-Powell terms to symmetrize the system or for instance by coupling equation \eqref{eqn.J} with a curl-cleaning field \cite{chiocchetti2021high}.


\section{Numerical scheme}
\label{sec:NumScheme}

\subsection{Mesh and notation}

\begin{figure}[!h]
	\centering
	\includegraphics[width=0.5\textwidth]{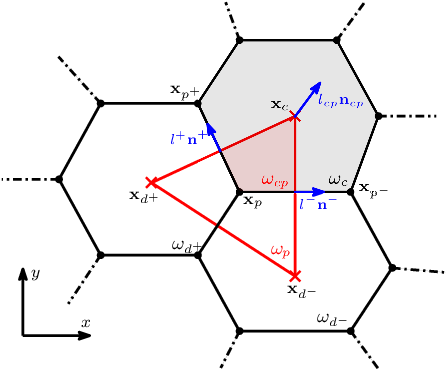}
	\caption{Mesh configuration: primal Voronoi cells (generators $\mathbf{x}_c$, volumes $|\omega_c|$) and dual Delaunay triangles (circumcenters $\mathbf{x}_p$, volumes $|\omega_p|$). The subcell $\omega_{pc}$ is highlighted in red.}
	\label{fig.mesh}
\end{figure}

In this work, we consider a two-dimensional computational domain denoted by $\Omega$ with boundary $\partial \Omega$ and defined by the position vector $\mathbf{x}=(x,y)$, partitioned into $\Nc$ non-overlapping Voronoi cells, see Fig.~\ref{fig.mesh}. 
Each cell $\omega_c$ is characterized by its generator $\x_c$, volume $|\omega_c|$. The grid is conforming in the sense that every interior edge in $\Omega \backslash \partial \Omega$ is shared by exactly two adjacent cells. Every vertex $\x_p$ of a Voronoi cell, represents the circumcenter of a triangular element $\omega_p$ of volume $|\omega_p|$, formed by the generators of the surrounding Voronoi cells. The resulting Delaunay triangulation is dual to the original Voronoi tessellation. The set of vertexes surrounding a Voronoi cell center $\x_c$ is denoted by $\mathcal{P}(c)$, while the set of cells whose centers form a triangle of circumcenter $\x_p$ is denoted by $\mathcal{C}(p)$. The set of faces $f$ which form the boundary $\partial \omega_c$ of a cell is $\mathcal{F}(c)$. The length of face $f$ of cell $\omega_c$ is labeled with $|\partial \omega_{fc}|$. The characteristic size of a Voronoi cell is measured by $h_c=\sqrt{|\omega_c|}$, and we take $h=\max_{\Omega} h_c$ as the characteristic mesh size over the entire computational domain. Furthermore, the subcell $\omega_{cp}$ is defined by connecting $\x_p$, $\x_c$, and the two midpoints of the edges $\x_{pp^-}$ and $\x_{pp^+}$. Finally, $l^{\pm}$ and $\n^{\pm}$ are the half-edge lengths and the normal vectors to the edges $\x_c \x_d^{\pm}$, respectively, so that the corner vector of $\omega_p$ originating from $\mathbf{x_c}$ is given by 
\begin{equation}
	l_{cp} \n_{cp} = l^- \n^- + l^+ \n^+, \quad \text{where} \quad l^{\pm} = \frac{1}{2} \norm{\x_c - \x_{d^{\pm}}}.
	\label{eqn.lpcnpc}
\end{equation}
Under these notations, one recovers a discrete version of the Gauss theorem such that  
\begin{equation}
	\sum_{c\in \mathcal{C}(p)} l_{cp} \n_{cp} = \mathbf{0}.
	\label{eqn.Gauss}
\end{equation}
Likewise, the tangent normal vector is given by
\begin{equation}
	l_{cp} \t_{cp} = ( -l_{cp} \n_{cp} \cdot \mathbf{e}_y, \, l_{cp} \n_{cp} \cdot \mathbf{e}_x ).
	\label{eqn.lpctpc}
\end{equation}
Obviously, it holds $l_{cp} \n_{cp} \cdot l_{cp} \t_{cp} = 0$. Notice that the outward pointing corner vectors of $\omega_c$ are simply the dual cell corner vectors $l_{cp} \n_{cp}$ with opposite sign for all $p \in \mathcal{P}(c)$, and thus the following identity holds:
\begin{equation}
	\sum_{p\in \mathcal{P}(c)} -l_{cp} \n_{cp} = \mathbf{0}.
	\label{eqn.GaussD}
\end{equation}
The same applies for the tangent vector $l_{cp} \t_{cp}$ in cell $\omega_c$.

\subsection{Compatible discrete curl-grad operators} \label{ssec.curlgrad}

For any arbitrary scalar field $\phi(\x)$ defined in the Voronoi cells $\omega_c$, one can define the gradient operator in the dual cells $\omega_p$ by 
\begin{equation}
	\mathbb{G}_p(\phi_c) := \nabla_h^c \phi_c = \frac{1}{\abs{\omega_p}} \sum_{c\in \mathcal{C}(p)} l_{cp} \n_{cp} \, \phi_c.
	\label{eqn.gradOp}
\end{equation}  
In contrast, for any arbitrary vector field $\boldsymbol{\psi}_p(\x)$ defined in the dual cells $\omega_p$, one can define the curl operator in the primal cells $\omega_c$ by
\begin{equation}
 	\mathbb{C}_c(\j_p) := \nabla_h^p \times \ppsi_p = -\frac{1}{\abs{\omega_c}} \sum_{p\in \mathcal{P}(c)} l_{cp} \t_{cp} \cdot \ppsi_p.
	\label{eqn.curlOp}
\end{equation}  

\begin{prop}[Compatible curl-grad operators] \label{prop1}
Under the definitions \eqref{eqn.gradOp} and \eqref{eqn.curlOp}, the following discrete curl-grad identity holds 
\begin{equation}
	\nabla_h^p \times \prn{\nabla_h^c \phi_c}_p = 0.
	\label{eqn.CurlGrad-P}
\end{equation}	
\end{prop}

\begin{proof}
	\begin{figure}[!h]
		\centering
		\includegraphics[width=0.4\textwidth]{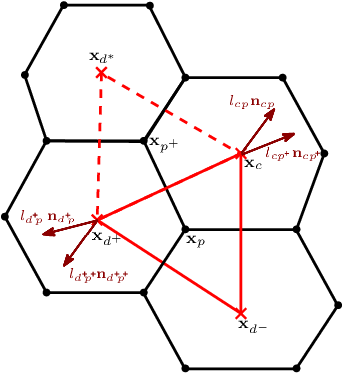}
		\caption{Local stencil around nodes $p$ and $p^+$, showing corner normals used in the proof of the discrete curl-grad identity~\eqref{eqn.CurlGrad-P}.}
		
		\label{fig.curlgrad}
	\end{figure}
	At the aid of Fig.~\ref{fig.curlgrad}, let us explicitly write the gradients obtained at nodes $p$ and $p^+$ using the operator \eqref{eqn.gradOp}:
	\begin{eqnarray}
		\mathbb{G}_p &=& \frac{1}{\abs{\omega_p}} \left( l_{cp} \n_{cp} \phi_c + l_{d^+p} \n_{d^+p} \phi_{d^+} + l_{d^-p} \n_{d^-p} \phi_{d^-} \right), \\
	   \mathbb{G}_{p^+} &=& \frac{1}{\abs{\omega_{p^+}}} \left( l_{cp^+} \n_{cp^+} \phi_c + l_{d^\star p^+} \n_{d^\star p^+} \phi_{d^\star} + l_{d^+ p^+} \n_{d^+ p^+} \phi_{d^+} \right). 
	\end{eqnarray}
    The above gradients provide two contributions to the curl operator of cell $\omega_c$, according to \eqref{eqn.curlOp}, that is
    \begin{equation}
    	-\frac{1}{\abs{\omega_c}} \left( \frac{1}{\abs{\omega_p}} \mathbb{G}_p \, l_{cp} \t_{cp} +\frac{1}{\abs{\omega_{p^+}}} \mathbb{G}_{p^+} \, l_{c p^+} \t_{c p^+} \right).
    \end{equation}
    The terms involving $\phi_c$ are
    \begin{equation}
    	-\frac{\phi_c}{\abs{\omega_c}} \left( \frac{1}{\abs{\omega_p}} \, l_{cp} \n_{cp} \cdot l_{cp} \t_{cp} +  \frac{1}{\abs{\omega_{p^+}}} \, l_{c {p^+}} \n_{c {p^+}} \cdot l_{c {p^+}} \t_{c {p^+}}  \right),
    \end{equation}
    which obviously vanish because of the orthogonal conditions $l_{cp} \n_{cp} \cdot l_{cp} \t_{cp}=0$ and $l_{c p^+} \n_{c p^+} \cdot l_{c p^+} \t_{c p^+}=0$, holding by construction \eqref{eqn.lpctpc}. Now, let us analyze the terms related to $\phi_{d^+}$:
    \begin{equation}
    	-\frac{\phi_{d^+}}{\abs{\omega_c}} \left( \frac{1}{\abs{\omega_p}} \, l_{d^+ p} \n_{d^+ p} \cdot l_{cp} \t_{cp} +  \frac{1}{\abs{\omega_{p^+}}} \, l_{d^+ {p^+}} \n_{d^+ {p^+}} \cdot l_{c {p^+}} \t_{c {p^+}}  \right),
    \end{equation}
    which is equivalent to
    \begin{equation}
    	-\frac{\phi_{d^+}}{\abs{\omega_c}} \left( \frac{1}{\abs{\omega_p}} \, l_{d^+ p} \n_{d^+ p} \times l_{cp} \n_{cp} +  \frac{1}{\abs{\omega_{p^+}}} \, l_{d^+ {p^+}} \n_{d^+ {p^+}} \times l_{c {p^+}} \n_{c {p^+}}  \right).
    	\label{eqn.phid}
    \end{equation}
    Using the following identities (see \cite{barsukow2025node} for the proof)
    \begin{equation}
    	l_{d^+ p} \n_{d^+ p} \times (-l_{cp} \n_{cp}) = -2 \, \abs{\omega_p},  \quad l_{d^+ {p^+}} \n_{d^+ {p^+}} \times (-l_{c {p^+}} \n_{c{p^+}}) = +2 \, \abs{\omega_{p^+}},
    \end{equation}
    the contribution \eqref{eqn.phid} reduces to zero. The same procedure applies for the remaining nodes of cell $\omega_c$. 
\end{proof}	

The consistency and accuracy of the discrete curl-grad operators is numerically verified by considering a computational domain $\Omega=[0;10]^2$ paved with $N_c=8090$ Voronoi cells with wall boundary conditions. The results are reported in Tab.~\ref{tab.curlgrad} for the scalar potential $\phi(\x)=\sin\left( \frac{2\pi x}{10}\right) \cos\left( \frac{2\pi y}{10}\right)$, showing first order of accuracy of the gradient operator ($\j_p=\nabla_h^c \phi_c$), and confirming the curl-free property.

\begin{table}[!h]
	\caption{Numerical validation of the discrete curl-grad operators on Voronoi cells.}
	\begin{center} 
		\begin{tabular}{c|cccc|c}
			$h$ & $L_{\infty}(j_1)$ & $\mathcal{O}(j_1)$ & $L_{\infty}(j_2)$ & $\mathcal{O}(j_2)$ & $|\nabla_h^p \times \j_p|_c$ \\
			\hline 
			2.771E-01 & 6.028E-02 &    - & 6.716E-02 &    - & 9.722E-15 \\
			1.847E-01 & 3.443E-02 & 1.38 & 3.515E-02 & 1.60 & 2.561E-14 \\
			1.386E-01 & 2.526E-02 & 1.08 & 2.599E-02 & 1.05 & 3.957E-14 \\
			1.111E-01 & 1.992E-02 & 1.08 & 2.084E-02 & 1.00 & 5.949E-14  
		\end{tabular}
	\end{center}
	\label{tab.curlgrad}
\end{table}

\begin{figure}[!h]
	\centering
	\includegraphics[width=0.7\textwidth]{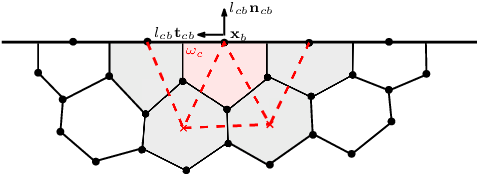}
	\caption{Portion of the computational grid around a node $b$ located on the boundary $\partial \Omega$.}
	\label{fig:boundary}
\end{figure}

\paragraph{Boundary conditions} Let $\j_b=(\nabla_h^c \phi_c)_b$ be the gradient defined on a boundary node $b$ of a cell $\omega_c$ with outward pointing corner normal vector $\n_{cb}$, and let $j_n=\j_b \cdot \n_{cb}$ be the prescribed boundary condition (see Fig.~\ref{fig:boundary}).

 The tangential component of $\j_b$, namely $j_t=\j_b \cdot \t_{cb}$, is a free degree of freedom that we choose to be given by
\begin{equation}
  j_t = -\frac{(\widetilde{\nabla_h^p \times \j_p})_c}{l_{cb}}, \qquad 	\widetilde{(\nabla_h^p \times \j_p)}_c = -\frac{1}{\abs{\omega_c}} \sum_{p\in \tilde{\mathcal{P}}(c)} l_{cp} \, \t_{cp} \cdot \j_p,
\end{equation}
where $\tilde{\mathcal{P}}(c)$ are the set of nodes of cell $\omega_c$ which do not lie on the boundary, and $\widetilde{(\nabla_h^p \times \j_p)}_c$ is the portion of the curl operator evaluated without accounting for the boundary contribution. It follows by construction that
\begin{equation}
	(\nabla_h^p \times \j_p) = \widetilde{(\nabla_h^p \times \j_p)}_c + l_{cb} \, \j_b \cdot  \t_{cb} = 0,
\end{equation}
so the curl-free property is preserved. Knowing $j_n$, prescribed by the boundary condition, and $j_t$, evaluated to ensure a zero curl operator, the vector $\j_b$ can be uniquely determined since $\n_{cb}$ and $\t_{cb}$ are known. The most common application is for wall boundary conditions where $j_n=0$, but other Dirichlet-type conditions can be easily imposed within the same framework.

\subsection{Discretization of the thermal impulse equation}
The semi-discrete evolution equation for the thermal impulse \eqref{eqn.J}, without source terms, writes on the dual triangular element $\omega_p$ as follows:
\begin{equation}
	\frac{\text{d} \j_p}{\text{d}t} = - \frac{1}{|\omega_p|} \sum \limits_{c \in \mathcal{C}(p)} l_{cp} \n_{cp} \varphi_c  - \frac{1}{|\omega_p|} \sum \limits_{c \in \mathcal{C}(p)} |\omega_{cp}| \, \mathbb{C}_c(\j_p) \times \vv_c.
	\label{eqn.jpsd}
\end{equation}
The thermal impulse is always initialized by means of a scalar potential $\phi_c$, that is $\j_p = \mathbb{G}_p(\phi_c)$. Consequently, by virtue of Property~\ref{prop1}, it holds $\mathbb{C}_c(\j_p) \times \vv_c = \mathbf{0}$. 

For the computation of the scalar quantity $\varphi_c=\j_c \cdot \vv_c + \temp_c$ we need a stable cell-centered value for $\j_c$. We propose to design a cell-solver, mimicking the nodal solvers widely adopted in Lagrangian schemes \cite{maire2009high,despres2005lagrangian,maire2011staggered,reale2011}, relying on the jump condition \cite{dhaouadi2024eulerian}
\begin{equation}
	\varphi_p-\varphi_c = z_p \, \left(\mathbf{\j_c}-\mathbf{\j_p}\right) \cdot \n_{pc}, \qquad z_p=\max(|\boldsymbol{\lambda}_p|)>0.
	\label{eqn.jump_varphi}
\end{equation}  
Conservation of the thermal impulse can be verified by summing the semi-discrete evolution equation \eqref{eqn.jpsd} over all the dual cells:
\begin{equation} 
	\frac{\text{d} }{\text{d}t} \left( \sum \limits_{p} |\omega_p| \, \j_p \right) = -\sum \limits_{p} \sum \limits_{c \in \mathcal{C}(p)} l_{cp} \n_{cp} \varphi_c 
	= -\sum \limits_{c} \sum \limits_{p \in \mathcal{P}(c)} l_{cp} \n_{cp} \varphi_c,
\end{equation}	
where we have switched the summation over cells and nodes since we are dealing with finite summations. A sufficient condition to ensure conservation is therefore
\begin{equation}
\sum \limits_{p \in \mathcal{P}(c)} l_{cp} \n_{cp} \varphi_c = 0.
	\label{eqn.consj}
\end{equation}
By employing the jump relation \eqref{eqn.jump_varphi} in the above expression, we obtain a vector equation for the unknown $\j_c$:
\begin{equation}
	\sum \limits_{p \in \mathcal{P}(c)} \mathbf{M}_{cp} \, \j_c = \sum \limits_{p \in \mathcal{P}(c)} \mathbf{M}_{cp} \, \j_p + l_{cp} \n_{cp} \varphi_p, \quad \mathbf{M}_{cp} = z_p \, \left( l_{cp}^- (\n_{cp}^- \otimes \n_{cp}^-) + l_{cp}^+ (\n_{cp}^+ \otimes \n_{cp}^+) \right),
	\label{eqn.jc}
\end{equation}
which is always solvable since the corner viscosity matrix $\mathbf{M}_{cp}$ is symmetric positive definite by construction \cite{maire_subcellForce}. Once $\j_c$ is known, the scalar $\varphi_c$ can be readily obtained to be used in \eqref{eqn.jpsd}.

\begin{prop}[Curl-free involution] \label{prop2}
Let $\mathbb{C}_p(\j_p)$ be a discrete curl operator defined in terms of \eqref{eqn.curlOp} which is acting on a vector $\j_p$ on the dual cell $\omega_p$:
\begin{equation}
	\mathbb{C}_p(\j_p) := \frac{1}{|\omega_p|} \, \sum_{c\in \mathcal{C}(p)} |\omega_{cp}| \, \mathbb{C}_c(\j_p).
	\label{eqn.curlp}
\end{equation}
The semi-discrete thermal impulse equation \eqref{eqn.jpsd} satisfies the involution \eqref{eqn.inv}, in the sense that
\begin{equation}
		\frac{\mathrm{d} }{\mathrm{d}t} \mathbb{C}_p(\j_p) = \mathbf{0}.
		\label{eqn.CFj}
\end{equation}

\end{prop}

\begin{proof}
Applying the curl operator \eqref{eqn.curlp} to the semi-discrete scheme \eqref{eqn.jpsd}, and recalling that $\mathbb{C}_c(\j_p) \times \vv_c = \mathbf{0}$, we get
\begin{eqnarray}
	\frac{\text{d} }{\text{d}t} \mathbb{C}_p(\j_p) &=& \frac{1}{|\omega_p|} \, \sum_{c\in \mathcal{C}(p)} |\omega_{cp}| \,   \mathbb{C}_c ( \, -\mathbb{G}_p(\varphi_c) \,)   = \mathbf{0},
\end{eqnarray}
where the flux contribution of $\varphi_c$ has been compactly rewritten using the gradient operator \eqref{eqn.gradOp}. Its curl is obviously zero by virtue of Property \ref{prop1}.
\end{proof}	

\begin{prop}[Local entropy inequality] The semi-discrete thermal impulse equation \eqref{eqn.jpsd}, with the scalar quantity $\varphi_c$ computed in terms of $\j_c$ obtained by the cell solver \eqref{eqn.jc}, yields a contribution to a local entropy inequality such that entropy increases.
\end{prop}

\begin{proof}
Looking at \eqref{eqn.pdq}, we see that the contribution of the semi-discrete thermal impulse equation \eqref{eqn.jpsd} to the semi-discrete entropy evolution can be computed from the dot product of \eqref{eqn.jpsd} with the conjugate variable $\beta_p :=(\partial_{\j} \rhos)_p=-\ar_p \j_p/\theta_p$, that is
	\begin{equation}
		\text{d}\rhos^{\, \j} = \beta \, \text{d}\j \quad \Rightarrow \quad |\omega_p| \,  \frac{\text{d} \rhos_p^{\, \j}}{\text{d}t} = -\beta_p \sum \limits_{c \in \mathcal{C}(p)} l_{cp} \n_{cp} (\varphi_c - \varphi_p),
		\label{eqn.sj}
	\end{equation} 
where $\rhos^{\, \j}$ denotes the contribution to entropy related to the thermal impulse equation only. The flux, which is written in fluctuation form thanks to the geometric identity \eqref{eqn.Gauss}, can be reformulated using the jump relation \eqref{eqn.jump_varphi} and the corner viscosity matrix $\mathbf{M}_{pc}$ defined in \eqref{eqn.jc}:
\begin{align}
	|\omega_p| \,  \frac{\text{d} \rhos_p^{\, \j}}{\text{d}t} &= \beta_p \sum \limits_{c \in \mathcal{C}(p)} \mathbf{M}_{cp} \, (\j_c-\j_p) \nonumber \\
	&= - \sum \limits_{c \in \mathcal{C}(p)} \underbrace{\frac{1}{2} (\beta_c - \beta_p) \, \mathbf{M}_{cp} \, (\j_c-\j_p)}_{\mathcal{S}_1} - \underbrace{\frac{1}{2} (\beta_c + \beta_p) \, \mathbf{M}_{cp} \, (\j_c-\j_p)}_{\mathcal{S}_2}.
	\label{eqn.sj2}
\end{align}
Because of the relation
\begin{equation}
	\int \limits_{\j_p}^{\j_c} \beta \, \mathbf{M} \, \text{d}\j = \int \limits_{\j_p}^{\j_c} \frac{\partial \rhos^{\j}}{\partial \j} \, \mathbf{M} \, \text{d}\j = \mathbf{M}_{cp} \, (\rhos^{\j}_c -\rhos^{\j}_p),
	\label{eqn.Roe_diss}
\end{equation}
the term $\mathcal{S}_2$ can be interpreted as the approximation of the jump of $s^{\j}$ by means of a trapezoidal rule for the path integral above, that is
\begin{equation}
	\int \limits_{\j_p}^{\j_c} \beta \, \mathbf{M} \, \text{d}\j \approx \frac{1}{2} (\beta_p + \beta_c) \, \mathbf{M}_{cp} \, (\j_c -\j_p),
	\label{eqn.trap_diss}
\end{equation}
which is a classical dissipation term of Laplace type, as usually adopted in finite volume schemes. The corner viscosity matrix does not depend on the path state, since it is computed at a frozen state (e.g. $z_{cp}=z_0$). For the term $\mathcal{S}_1$, the jump in the conjugate variables is now rewritten as a jump in the conservative variables through the Hessian matrix of the entropy potential which satisfies the Roe property
\begin{equation}
	\boldsymbol{\mathcal{H}}_{cp}^{\, \j} \cdot (\j_c-\j_p) = \beta_c-\beta_p.
	\label{eqn.HRoe}
\end{equation}
The Hessian matrix of the energy potential \eqref{eqn.rhoE}, restricted to the contribution of the thermal impulse, is analytically given by
\begin{equation}
\partial^2_{\j \j} \left(\frac{1}{2} \ar \| \j \|^2 \right)_{cp} = \left( \begin{array}{cc}
		\ar_{cp} & 0 \\ 0 & \ar_{cp}
	\end{array} \right),
\end{equation} 
which is symmetric and positive definite since $\ar_{cp}>0$. Because of the convexity of the energy potential \eqref{eqn.rhoE}, the Hessian matrix of the entropy potential $\boldsymbol{\mathcal{H}}_{cp}^{\, \j}$ is symmetric and negative definite. Using the relations \eqref{eqn.HRoe} and \eqref{eqn.trap_diss} in \eqref{eqn.sj2} yields
\begin{equation}
	|\omega_p| \,  \frac{\text{d} \rhos_p^{\, \j}}{\text{d}t} =\underbrace{-\sum \limits_{p \in \mathcal{P}(c)} \frac{1}{2} \, (\j_c-\j_p)^\top \, \boldsymbol{\mathcal{H}}_{cp}^{\, \j} \, \mathbf{M}_{cp} \, (\j_c-\j_p)}_{\geq 0} + \underbrace{\sum \limits_{p \in \mathcal{P}(c)} \mathbf{M}_{cp} \, (\rhos^{\j}_c -\rhos^{\j}_p)}_{\geq 0}.
	\label{eqn.sj3}
\end{equation}	
We can then conclude that \eqref{eqn.sj} leads to
\begin{equation}
	|\omega_p| \, \frac{\text{d} \rhos_p^{\, \j}}{\text{d}t} \geq 0,
\end{equation}
showing that the semi-discrete thermal impulse equation \eqref{eqn.jpsd} increases the corresponding contribution $\rhos_p^{\, \j}$ to local entropy within the dual cell $\omega_p$.
\end{proof}

By accounting for the dissipation source terms, the fully discrete implicit-explicit scheme for the evolution of $\j$ is given by
\begin{equation}
	\j^{n+1}_p = \j^{n}_p - \frac{\dt}{|\omega_p|} \sum \limits_{c \in \mathcal{C}(p)} l_{pc} \n_{pc} \varphi_c^n  - \frac{\dt}{|\omega_p|} \sum \limits_{c \in \mathcal{C}(p)} |\omega_{pc}| \, \mathbb{C}_c(\j^{n}_p) \times \vv_c^n - \dt \frac{\j_p^{n+1}}{\tau}.
	\label{eqn.jpfd}
\end{equation}

\begin{prop}[Asymptotic-preservation of the Fourier limit] \label{prop3}
Let us introduce the $k$-th order Chapman-Enskog expansion of the thermal impulse $\j$ in powers of the stiffness parameter $\tau$:
\begin{equation}
	\j = \j_{(0)} + \tau  \j_{(1)} + \tau^2 \j_{(2)} + \ldots + \mathcal{O}(\tau^k).
	\label{eqn.CE}
\end{equation} 
Assuming well-prepared initial data, i.e. $\j^{n}_{p,(0)}=\mathcal{O}(\tau^2)$, the fully discrete scheme \eqref{eqn.jpfd} is asymptotic preserving in the sense that in the stiff relaxation limit for $\tau \to 0$ the heat flux \eqref{eqn.heatflux} is retrieved also at the fully discrete level.
\end{prop}

\begin{proof}
Application of the asymptotic expansion \eqref{eqn.CE} to the fully discrete scheme \eqref{eqn.jpfd} and collection of the like powers of $\tau$ up to leading order terms yield
\begin{align}
\j^{n+1}_{p,(0)} = \j^{n}_{p,(0)} &- \frac{\dt}{|\omega_p|} \sum \limits_{c \in \mathcal{C}(p)} l_{pc} \n_{pc} (\j_{c,(0)}^n \cdot \vv_c^n+\temp_c^n) \nonumber \\
& - \frac{\dt}{|\omega_p|} \sum \limits_{c \in \mathcal{C}(p)} |\omega_{pc}| \, \mathbb{C}_c(\j^{n}_{p,(0)}) \times \vv_c^n - \dt \frac{\j_{p,(0)}^{n+1}}{\tau} - \dt \, \j_{p,(1)}^{n+1} + \mathcal{O}(\tau^2).	
\label{eqn.Jexpansion}
\end{align}
From the leading order term $\tau^{-1}$ it follows
\begin{equation}
	\j^{n+1}_{p,(0)} = \mathcal{O}(\tau^2).
	\label{eqn.j0}
\end{equation}
Since the computation of $\j_c$ using the solver \eqref{eqn.jc} does not depend on $\tau$ and $\j^{n}_{p,(0)}=\mathcal{O}(\tau^2)$, we can conclude that also $\j^{n}_{c,(0)}=\mathcal{O}(\tau^2)$. Using this fact and the above result \eqref{eqn.j0} in the scheme \eqref{eqn.Jexpansion}, the terms of order $\tau^0$ are
\begin{equation}
	\j_{p,(1)}^{n+1} = - \frac{1}{|\omega_p|} \sum \limits_{c \in \mathcal{C}(p)} l_{pc} \n_{pc} \temp_c^n.
	\label{eqn.j1}
\end{equation}
We insert \eqref{eqn.j0} and \eqref{eqn.j1} in the expansion \eqref{eqn.CE}, hence obtaining
\begin{equation}
	\j^{n+1}_{p} = - \frac{\tau}{|\omega_p|} \sum \limits_{c \in \mathcal{C}(p)} l_{pc} \n_{pc} \temp_c^n.
\end{equation}
Consequently, in the stiff relaxation limit $\tau \to 0$, the heat flux $\mathbf{q}=\ar \, \temp \, \j$ in the energy equation \eqref{eqn.en} becomes
\begin{equation}
	\mathbf{q}_p^{n} = - \tau \, \ar_p \, \temp_p \frac{1}{|\omega_p|} \sum \limits_{c \in \mathcal{C}(p)} l_{pc} \n_{pc} \temp_c^n,
\end{equation}
which is a consistent discretization of \eqref{eqn.heatflux}.
\end{proof}	

\subsection{Discretization of mass, momentum and energy equations}
Because of the staggered discretization of the thermal impulse equation, the fluxes $\h=\f+\g$ of the remaining equations of system \eqref{eqn.heat_system} are split as follows:
\begin{equation}
\f = \left( \begin{array}{l}
		\f^{\rho} \\ \f^{\rho \vv} \\ \f^{E}
	\end{array} \right) = \left( \begin{array}{c}
		\rho \vv \\ \rho \vv \otimes \vv + p \mathbf{I} \\ (E+p) \, \vv
	\end{array} \right), \quad \g=\left( \begin{array}{l}
	\g^{\rho} \\ \g^{\rho \vv} \\ \g^{E}
\end{array} \right) = \left( \begin{array}{c}
0 \\ -\ar \| \j \|^2 \, \mathbf{I} + \ar \, \j \otimes \j \\ (-\ar \| \j \|^2 \, \mathbf{I} + \ar \, \j \otimes \j) \, \vv
\end{array} \right).
	\label{eqn.flux}
\end{equation}  
The semi-discrete finite volume scheme on Voronoi grids is then given by
\begin{subequations}
	\label{eqn.fv}
	\begin{align}
		\frac{\text{d} \rho_c}{\text{d}t} &= - \frac{1}{|\omega_c|} \sum \limits_{f \in \mathcal{F}(c)} |\partial \omega_{fc}| \, \hat{\f}_{f}^{\rho} \cdot \n_{fc}, \label{eqn.rho_fv}  \\
		\frac{\text{d} (\rho \vv)_c}{\text{d}t} &= - \frac{1}{|\omega_c|} \sum \limits_{f \in \mathcal{F}(c)} |\partial \omega_{fc}| \,  \left( \hat{\f}_{f}^{\rho \vv} + \frac{1}{2}\sum \limits_{p \in \mathcal{P}(f)} {\g}_{pf}^{\rho \vv} \right) \cdot \n_{fc}, \label{eqn.mom_fv} \\
		\frac{\text{d} E_c}{\text{d}t} &= - \frac{1}{|\omega_c|} \sum \limits_{f \in \mathcal{F}(c)} |\partial \omega_{fc}| \,  \left( \hat{\f}_{f}^{E} + \frac{1}{2}\sum \limits_{p \in \mathcal{P}(f)} {\g}_{pf}^{E} \right) \cdot \n_{fc}, \label{eqn.en_fv}
	\end{align}
\end{subequations}
where we use a Rusanov--type numerical flux \cite{Rusanov:1961a}, which, for a generic variable $q(\x,t)$, writes 
\begin{equation}
	\hat{\f}_{f}^{q} \cdot \n_{fc} = \frac{1}{2} \left( \f^{q}_c + \f^{q}_d \right) \cdot \n_{f} - \frac{1}{2} z_{f} (q_d-q_c),
	\label{eqn.rusanov}
\end{equation}
being $\omega_d$ the cell which shares the face $f$ with cell $\omega_c$. The numerical dissipation coefficient is given by $z_{f} = \max\left(  \max(|\boldsymbol{\lambda}_c|), \, \max(|\boldsymbol{\lambda}_d|) \right)$. 

\subsection{Thermodynamic compatibility}
The semi-discrete scheme \eqref{eqn.fv} is compliant with the First Law of Thermodynamics, since total energy is conserved by \eqref{eqn.en_fv}. To be compliant also with the Second Law of Thermodynamics, we need to recover the entropy balance \eqref{eqn.eta} by means of a discrete analogue of \eqref{eqn.pdq}:
\begin{equation}
	\begin{array}{c}
		\displaystyle \temp_c \cdot \text{d}\rhos_c = - r_c \cdot \text{d}\rho_c - \vv_c \cdot \text{d}(\rho \vv)_c - b_c \cdot \text{d}\j_c + 1 \cdot \text{d}E_c, \\ [2mm]
		\displaystyle \temp_c= (\partial_\rhos E)_c, \quad r_c = (\partial_{\rho}E)_c, \quad \vv_c = (\partial_{\rho \mathbf{u}}E)_c, \quad b_c = (\partial_{\j}E)_c.
	\end{array}
		\label{eqn.pdq_c}
\end{equation}
For the moment, let us assume to have a semi-discrete evolution equation for the cell quantity $\j_c$ with associated numerical fluxes $\hat{\h}^{\j}$. Following the general framework of \cite{Abgrall2018}, in our scheme \eqref{eqn.fv} we introduce a conservative
modification of the numerical fluxes across the face $f$ shared between cells $\omega_c$ and $\omega_d$:
\begin{equation}
	\tilde{\h}_f \cdot \n_{fc} = \hat{\h}_f \cdot \n_{fc} - \as ( \p_d - \p_c ), \qquad \hat{\h}_f:=\hat{\f}_f + \frac{1}{2}\sum \limits_{p \in \mathcal{P}(f)} {\g}_{pf},
	\label{eqn.fluxS}
\end{equation}
with $\as$ a scalar face correction factor that multiplies the jump in the thermodynamic dual variables $\p=(-r/\theta, \, -\mathbf{u}/\theta, \, -b/\theta, \, 1/\theta)^\top$. To determine $\as$, we rely on conservation principles \cite{HTCAbgrall}, hence requiring that the sum of all the fluctuations across the element interface $f$ is equal to the flux difference of the entropy equation which we want to satisfy:
\begin{equation}
	\p_c \cdot \left( \tilde{\h}_f \cdot \n_{fc} - \h_c \cdot \n_{fc} \right) + \p_d \cdot \left( \h_d \cdot \n_{fc} - \tilde{\h}_f \cdot \n_{fc} \right) = \left(  \f^{\rhos}_d -  \f^{\rhos}_c \right) \cdot \n_{fc},
	\label{eqn.cons}
\end{equation}
where $\f^{\rhos}=\rhos \vv+\ar\,\j $ and we recall that $\h=\f+\g$. Inserting the flux definition \eqref{eqn.fluxS} in the conservation condition \eqref{eqn.cons} yields the sought correction factor, that is
\begin{equation}
	\as = \frac{\left(  \f^{\rhos}_d -  \f^{\rhos}_c \right) \cdot \n_{fc} + \left(\hat{\h}_f \cdot \n_{fc}\right) \cdot (\p_d-\p_c) - (\p_d \cdot \h_d - \p_c \cdot \h_c) \cdot \n_{fc} }{(\p_d-\p_c)^2}.
	\label{eqn.alpha}
\end{equation}
In the case of vanishing denominator in \eqref{eqn.alpha}, we set $\as=0$. 
We still have to determine how to obtain the fluxes $\hat{\h}^{\j}$ for the evolution equation of $\j_c$.

\subsection{Virtual semi-discrete evolution of the thermal impulse equation} \label{ssec.jc}
\begin{figure}[!h]
	\centering
	\includegraphics[width=0.6\textwidth]{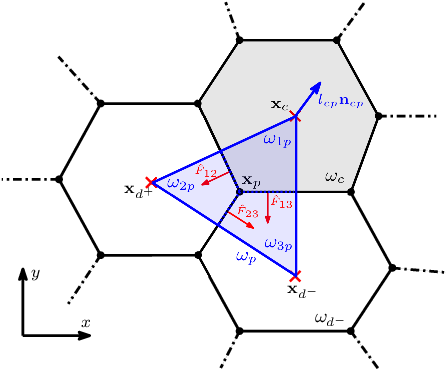}
	\caption{Subgrid of the dual cell $\omega_p$ with internal subfluxes $\hat{F}^\j$.}
	\label{fig.subgrid}
\end{figure}

Let us introduce a subgrid for each dual triangular volume $\omega_p$, as depicted in Fig.~\ref{fig.subgrid}. Because of the Voronoi tessellation, each triangle has always $N_m=3$ triangular subcells $\omega_{mp}$, $N_f=3$ internal subfaces and $N_{f^b}=3$ boundary subfaces (i.e. the sum of the two subfaces which lie on $\partial \omega_p$ for each subcell). For each subcell, $\mathcal{F}(m)$ and $\mathcal{F}^b(m)$ are the set of internal and boundary subfaces, respectively. 
According to the numbering shown in Fig.~\ref{fig.subgrid}, the adjacency matrix $\A_{[N_m \times N_f]}$ writes
\begin{equation}
	\A = \left( \begin{array}{ccc}
		\phantom{+}1 & \phantom{+}0 & \phantom{+}1 \\
		-1 & \phantom{+}1 & \phantom{+}0 \\
		\phantom{+}0 & -1 & -1
	\end{array} \right), \qquad \text{with} \qquad A_{mr} = \left\{ \begin{array}{rl}
	1 & \text{if} \quad m < r \\
	-1 & \text{if} \quad m > r \\
	0 & \text{otherwise}
\end{array} \right.,
\label{eqn.Amat}
\end{equation}
where $m$ and $r$ are the indexes of the subcells sharing the subface $f$. To obtain a semi-discrete evolution equation for $\j_c$, we want to reformulate the update \eqref{eqn.jpsd} as a subcell finite volume scheme of the form
\begin{equation}
	\frac{\text{d} \j_m}{\text{d}t} = -\frac{1}{|\omega_m|} \left( \sum \limits_{f \in \mathcal{F}(m)}  \hat{F}_{fm}^\j + \sum \limits_{f \in \mathcal{F}^b(m)} \hat{B}_{fm}^\j \right), \qquad m \in \{1,N_m\},
	\label{eqn.scfv}
\end{equation}
where $\hat{F}_{fm}^\j$ and $\hat{B}_{fm}^\j$ are the set of internal and boundary fluxes, respectively. The above subcell finite volume scheme \eqref{eqn.scfv} can be written in matrix notation \cite{Vilar2024} as
\begin{equation}
	\begin{array}{c}
		\displaystyle \frac{\text{d} \J_p}{\text{d}t} = -\D_p^{-1} \, \left( \A \hat{\F}_p + \hat{\B}_p \right), \\ [2mm]
		\J_p=\{\j_m\}, \, \,  \D_p=\{|\omega_m|\}, \, \, \hat{\F}_p=\{\hat{F}_{f}^\j\}, \, \, \hat{\B}_p=\{\hat{B}_{f^b}^\j\}.
	\end{array} 
\label{eqn.scfv2}
\end{equation}
We can notice that the boundary fluxes are known and, indeed, they write
\begin{equation}
	\hat{\B}_p=\{l_{pc} \n_{pc} \varphi_c\}, \qquad c \in \mathcal{C}(p).
\end{equation}
The unknown internal subfluxes $\hat{\F}_p$ can be determined by imposing the equivalence of the subcell scheme \eqref{eqn.scfv2} with the semi-discrete evolution of $\j_p$ given by \eqref{eqn.jpsd}:
\begin{equation}
	-\D_p^{-1} \, \left( \A \hat{\F}_p + \hat{\B}_p \right) = \boldsymbol{\Phi}_p,
	\label{eqn.scfluxA}
\end{equation}
where the abbreviation $\boldsymbol{\Phi}_p$ refers to the right hand side of \eqref{eqn.jpsd}, also called residual. As pointed out in \cite{Abgrall_subflux2018}, the adjacency matrix is singular, therefore we proceed by computing the graph Laplacian $\L_{[N_m \times N_m]}$ and its pseudo-inverse as
\begin{equation}
	\L=\A \A^{\top}, \qquad \L^{-1} = (\L + \nu \Pi)^{-1} - \frac{\Pi}{\nu},
	\label{eqn.GL}
\end{equation}
with $\{\nu \in \mathbb{R} \, | \, \nu \neq 0\}$ and $\Pi=(\mathds{1} \otimes \mathds{1})/N_m$. Eventually, we can solve the relation \eqref{eqn.scfluxA} and compute the subcell fluxes by means of \eqref{eqn.GL}, that is
\begin{equation}
	\hat{\F}_p = -\A^{\top} \L^{-1} \left( \D_p \Phi_p + \hat{\B}_p \right).
	\label{eqn.scflux} 
\end{equation}
Now, these subfluxes can be readily used to write a virtual semi-discrete evolution for the thermal impulse equation in the cell:
\begin{equation}
	\frac{\text{d} \j_c}{\text{d} t}  = -\frac{1}{|\omega_c|} \sum \limits_{p \in \mathcal{P}(c)} \boldsymbol{\sigma}_{pc}  \hat{\mathbf{F}}_{pc} + \frac{\j_c}{\tau},
	\label{eqn.sdjc}
\end{equation}
where $\hat{\mathbf{F}}_{pc} \in \hat{\mathbf{F}}_{p}$ is the set of subfluxes of cell $\omega_p$ that share node $p$ with cell $c$. $\boldsymbol{\sigma}_{pc}$ is a sign function such that makes the subfluxes projected along the outward pointing normal $\n_{pc}$. For instance, in Fig.~\ref{fig.subgrid}, for cell $\omega_{c}$ and point $p$, we have $\hat{\mathbf{F}}_{pc} = (\hat{F}_{12}, \hat{F}_{13})$ and $\boldsymbol{\sigma}_{pc}=(1,1)$, which corresponds to entries $A_{11}$ and $A_{13}$ of the adjacency matrix \eqref{eqn.Amat}. By defining $\hat{\h}^{\j}=\boldsymbol{\sigma}_{pc}  \hat{\mathbf{F}}_{pc}$, the virtual semi-discrete evolution equation \eqref{eqn.sdjc} is suitable to comply with the thermodynamically compatible correction \eqref{eqn.alpha}.

\begin{prop}[Thermodynamically compatible semi-discrete scheme] Assuming impermeable boundary conditions $\int \limits_{\partial \Omega} \h \cdot \n \, \text{dS} =\mathbf{0}$ and smooth flows, the semi-discrete scheme \eqref{eqn.fv} with the modified fluxes \eqref{eqn.fluxS} and the virtual semi-discrete evolution equation \eqref{eqn.sdjc} satisfies the semi-discrete entropy equation
\begin{equation}
	\frac{\text{d}\rhos_c}{\text{d}t} +  \frac{1}{|\omega_c|} \sum \limits_{f \in \mathcal{F}(c)} |\partial \omega_{fc}| \, \mathcal{F}_{f}^{\rhos} \cdot \n_{fc} = \frac{\ar_c}{ \temp_c\,\tau}\norm{\j_c}^2, \label{eqn.sdeta}
\end{equation}
with the compatible fluxes given by
\begin{equation}
	\mathcal{F}_{f}^{\rhos} \cdot \n_{fc} = \frac{1}{2} \left(  \, (\p_c \cdot \tilde{\h}_f + \f_c^{\rhos} - \p_c \cdot \h_c) + (\p_d \cdot \tilde{\h}_f + \f_d^{\rhos} - \p_d \cdot \h_d) \, \right) \cdot \n_{fc}. 
	\label{eqn.compEntFlux}
\end{equation}
\end{prop}	

\begin{proof}
The semi-discrete scheme \eqref{eqn.fv} compactly writes
\begin{equation}
	\frac{\text{d}\w_c}{\text{d}t} + \frac{1}{|\omega_c|} \sum \limits_{f \in \mathcal{F}(c)} |\partial \omega_{fc}| \, \tilde{\h}_{f} \cdot \n_{fc},
\end{equation}
with $\w_c=(\rho_c, \rho \vv, E, \j)_c^\top$ and the numerical fluxes $\tilde{\h}_{f}$ which follows from the definition \eqref{eqn.fluxS}. After dot multiplying the above equation by the dual variables $\p_c$, and adding and subtracting the term $\frac{1}{2} \p_d \cdot (\tilde{\h}_f \cdot \n_{fc})$, we get
\begin{equation}
	\frac{\text{d}\rhos_c}{\text{d}t} + \frac{1}{2 \, |\omega_c|} \sum \limits_{f \in \mathcal{F}(c)} |\partial \omega_{fc}| \, (\p_c+\p_d) \cdot \tilde{\h}_{f} \cdot \n_{fc} + (\p_c-\p_d) \cdot \tilde{\h}_{f} \cdot \n_{fc} =\frac{\ar \| \j_c\|^2}{\theta_c \, \tau}.
\end{equation}
The source term is automatically compatible, since it is the pointwise discrete analogue of the continuous model. Let us use the condition \eqref{eqn.cons} to rewrite the term $(\p_c-\p_d) \cdot \tilde{\h}_f \cdot \n_{fc}$ as
\begin{align}
	\frac{\text{d}\rhos_c}{\text{d}t} & + \frac{1}{2 \, |\omega_c|} \sum \limits_{f \in \mathcal{F}(c)} |\partial \omega_{fc}| \, \left( (\p_c+\p_d) \cdot \tilde{\h}_{f} + (\f_d^{\rhos}+\f_c^{\rhos}) + (\p_c \h_c - \p_d \h_d) \right) \cdot \n_{fc} \nonumber\\
	&= \frac{\ar \| \j_c\|^2}{\theta_c \, \tau}.
\end{align}
We now add, to the left hand side of the above entropy equation, the term 
\begin{equation}
\frac{1}{2} \sum \limits_{f \in \mathcal{F}(c)} |\partial \omega_{fc}| (\f_c^{\rhos}-\p_c \h_c) \cdot \n_{fc}=0,
\end{equation}
which is zero according to the discrete Gauss theorem, hence yielding
\begin{align}
	\frac{\text{d}\rhos_c}{\text{d}t} & + \frac{1}{2 \, |\omega_c|} \sum \limits_{f \in \mathcal{F}(c)} |\partial \omega_{fc}| \,  (\p_c+\p_d) \cdot \tilde{\h}_{f} \cdot \n_{fc} \nonumber \\
	& + \frac{1}{2 \, |\omega_c|} \sum \limits_{f \in \mathcal{F}(c)} |\partial \omega_{fc}| \, (\f_d^{\rhos}+\f_c^{\rhos}) \cdot \n_{fc} \nonumber \\
	&- \frac{1}{2 \, |\omega_c|} \sum \limits_{f \in \mathcal{F}(c)} |\partial \omega_{fc}| \,  (\p_c \h_c + \p_d \h_d)  \cdot \n_{fc} \nonumber\\
	&= \frac{\ar \| \j_c\|^2}{\theta_c \, \tau},
\end{align}
where we can recognize the flux $\mathcal{F}_f \cdot \n_{fc}$ given by \eqref{eqn.compEntFlux}.
\end{proof}	

\section{Numerical results}
\label{sec:Results}
We present a suite of test cases to numerically show the structure-preserving properties of our new finite volume scheme, which is labeled with SPFV. If not otherwise stated, we set $\gamma=2$ and $c_v=1$. The thermal impulse $\j$ is always initialized by taking the discrete gradient of a thermal potential $\phi$, according to its definition \cite{dhaouadi2024eulerian}. Slip wall boundary conditions are imposed by default as presented in Sec.~\ref{ssec.curlgrad}, and the initial condition is given in terms of the primitive variables $\mathbf{v}=(\rho,u_1,u_2,p,j_1,j_2)$. All the numerical results presented in this section are done for $\ar=\kappa^2/\rho$ with $\kappa=1$ set by default.

\subsection{Convergence analysis} \label{ssec.conv}
To verify the accuracy of the novel structure preserving finite volume scheme, we set up a manufactured solution for the model \eqref{eqn.heat_system} in the absence of source terms. The computational domain is the square $\Omega=[0,10]^2$, and the fluid is assigned a unit density and pressure with the following velocity field and thermal potential:
\begin{equation}
	\label{eqn.ICconv}
	\left(\begin{array}{c}
			\displaystyle u_1 \\ 	\displaystyle u_2
	\end{array}\right) = \left( \begin{array}{c}
	\displaystyle -\frac{\beta e^{\frac{1}{2}(1-r^2)}}{2\pi}y \\ 	\displaystyle +\frac{\beta e^{\frac{1}{2}(1-r^2)}}{2\pi}x
\end{array} \right), \qquad \phi = -\frac{\lambda e^{\frac{1}{2}(1-r^2)}}{2\pi}.
\end{equation}
By choosing $\beta=1$ and $\lambda=1$, the governing equations are balanced by adding only one artificial source term in the energy equation:
\begin{equation}
	\label{eqn.SEconv}
	S(\rho E) = -\frac{e^{\frac{1}{2}(1-r^2)}(r^2-2)}{2\pi}. 
\end{equation}
The simulation is run on a sequence of successively refined Voronoi meshes of size $h$ until time $t_f=0.1$. The results are reported in Tab.~\ref{tab.conv} exhibiting first order of accuracy, as expected. The errors related to the curl of $\j$ at the final time of each simulation are shown as well, serving as  numerical evidence to the curl-free property of the novel scheme.

\begin{table}[H]
	\begin{center}
		{
			\caption{Experimental order of convergence for the manufactured solution presented in Sec.~\ref{ssec.conv}. The errors are measured in the $L_2$ norm and refer to the variables $\rho$ (density), $u_1$ (horizontal velocity), $p$ (pressure) and $j_1$ (horizontal thermal impulse) at time $t_f=0.1$. The $L_{\infty}$ norm of the curl of $\j$ is also reported.} \label{tab.conv}
			\begin{small}
				\renewcommand{\arraystretch}{1.1}
				\begin{tabular}{c|cccccc|c}
					$h$ & ${L_2}(\rho)$ & $\mathcal{O}(\rho)$ & ${L_2}(u)$ & $\mathcal{O}(u)$ & ${L_2}(j_1)$ & $\mathcal{O}(j_1)$ & ${L_{\infty}(\nabla \times \j)_z}$ \\
					\hline
					3.75E-01 & 1.05E-03 &    - & 5.48E-02 &    - &  4.99E-02 &    - & 7.22E-16 \\
					2.95E-01 & 5.94E-04 & 1.45 & 3.58E-02 & 1.09 &  3.16E-02 & 1.17 & 2.43E-15 \\
					2.25E-01 & 4.21E-04 & 1.24 & 2.71E-02 & 1.00 &  2.35E-02 & 1.05 & 4.00E-15 \\
					1.84E-01 & 3.33E-04 & 1.08 & 2.19E-02 & 0.99 &  1.89E-02 & 1.02 & 7.27E-15 \\
				\end{tabular}
			\end{small}
		}
	\end{center}
\end{table}

\subsection{Riemann problems} \label{ssec.RP}
\begin{figure}[!htbp]
	\begin{center}
		\begin{tabular}{cc}
			\includegraphics[trim= 5 5 5 5,clip,width=0.45\textwidth,keepaspectratio=true]{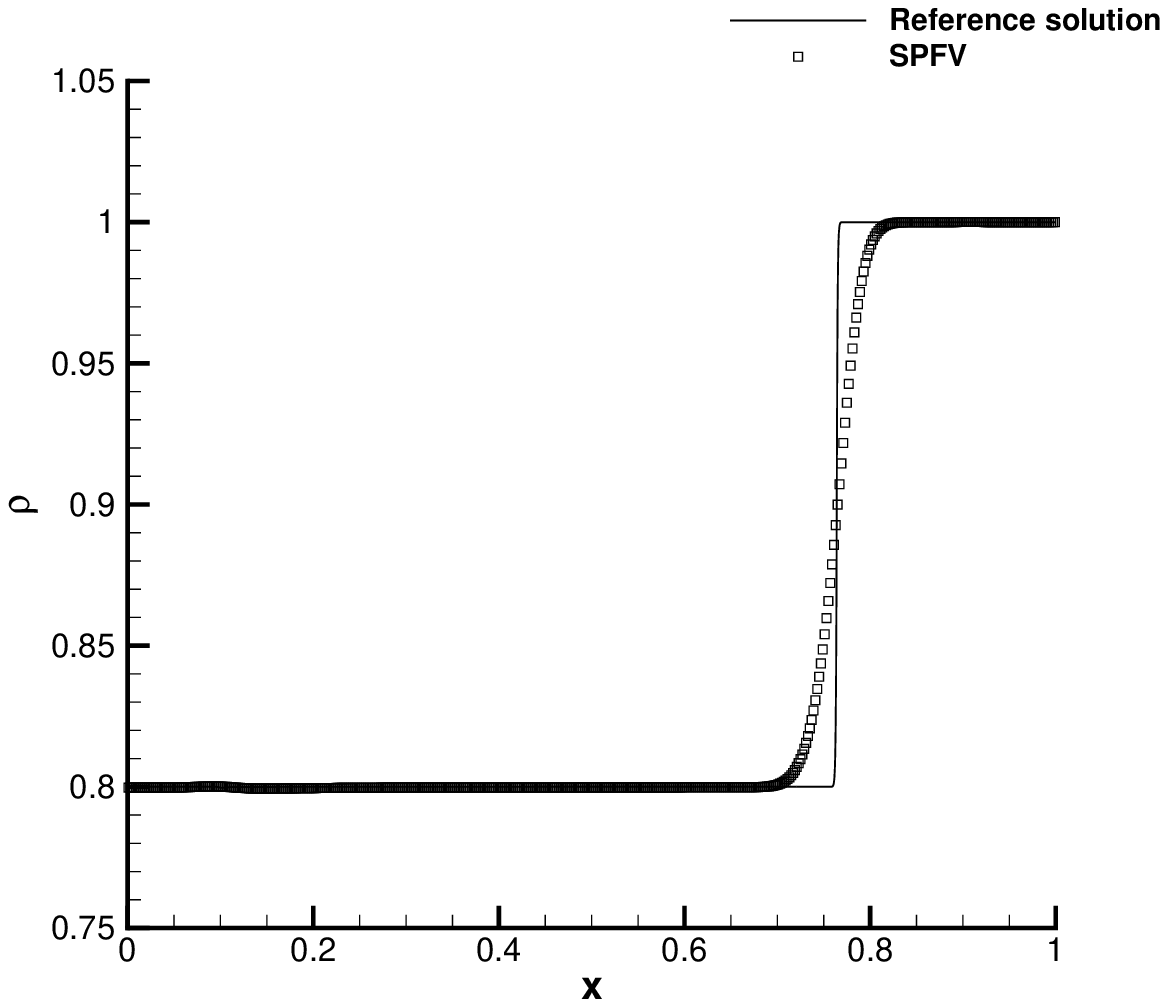} &
			\includegraphics[trim= 5 5 5 5,clip,width=0.45\textwidth,keepaspectratio=true]{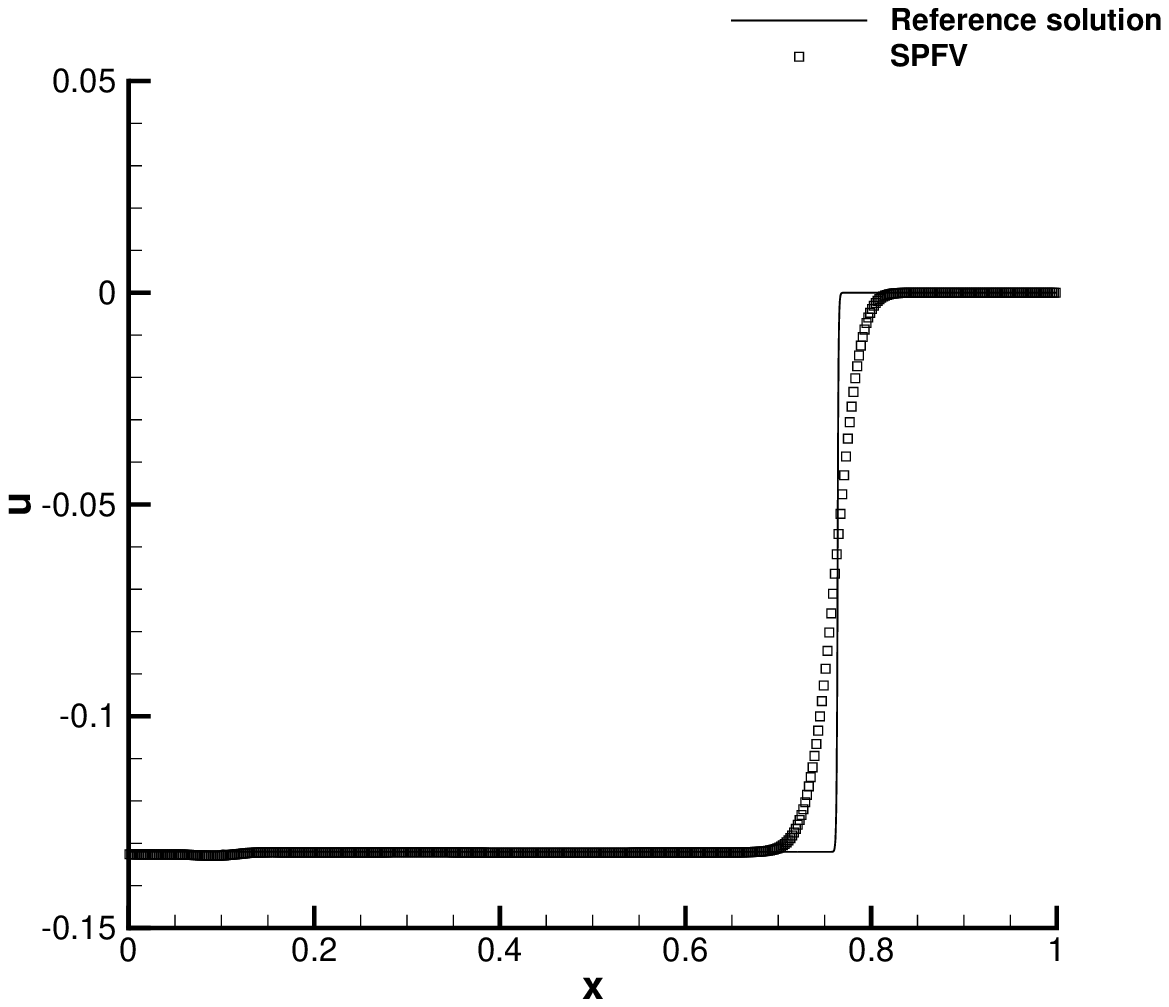} \\
			\includegraphics[trim= 5 5 5 5,clip,width=0.45\textwidth,keepaspectratio=true]{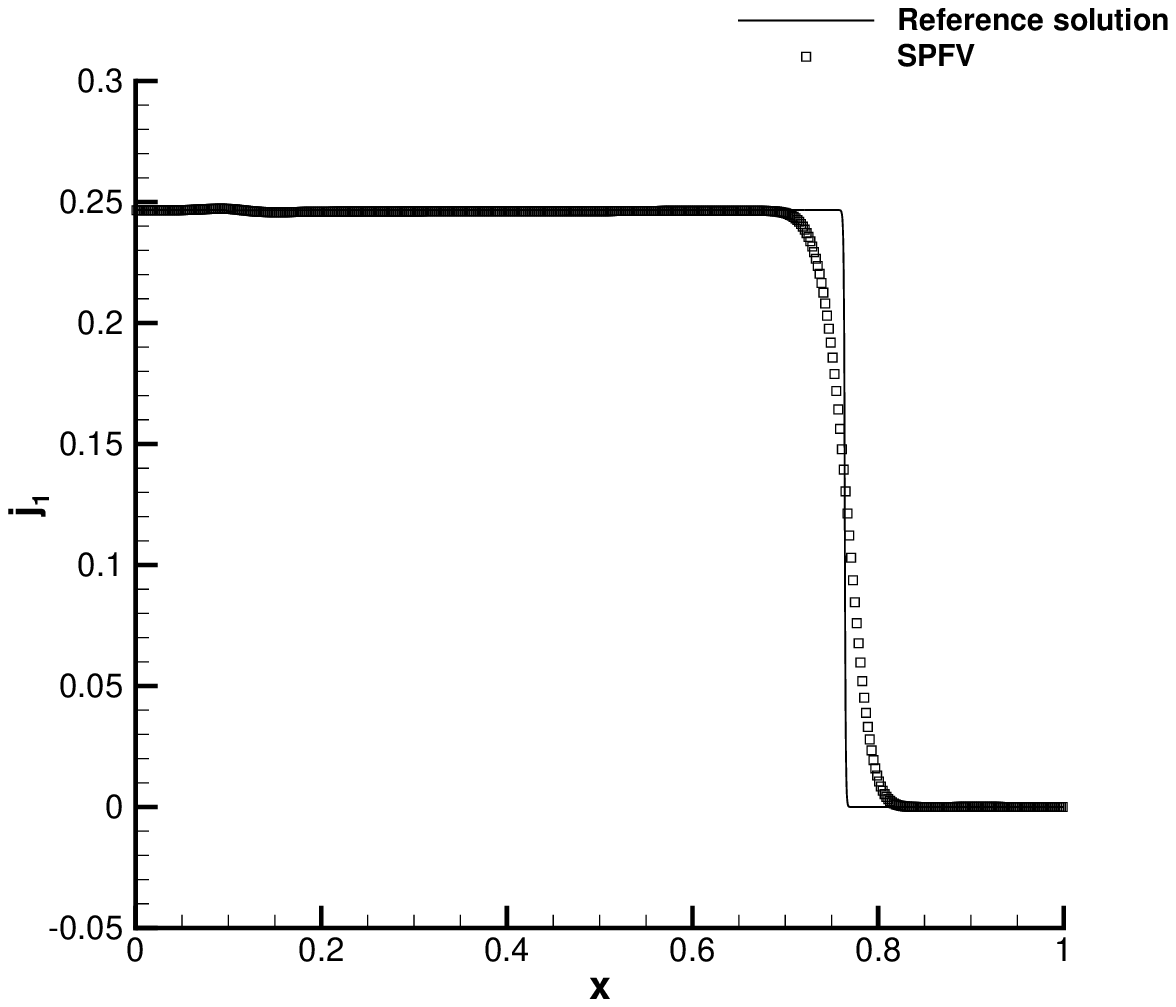} &
			\includegraphics[trim= 5 5 5 5,clip,width=0.45\textwidth,keepaspectratio=true]{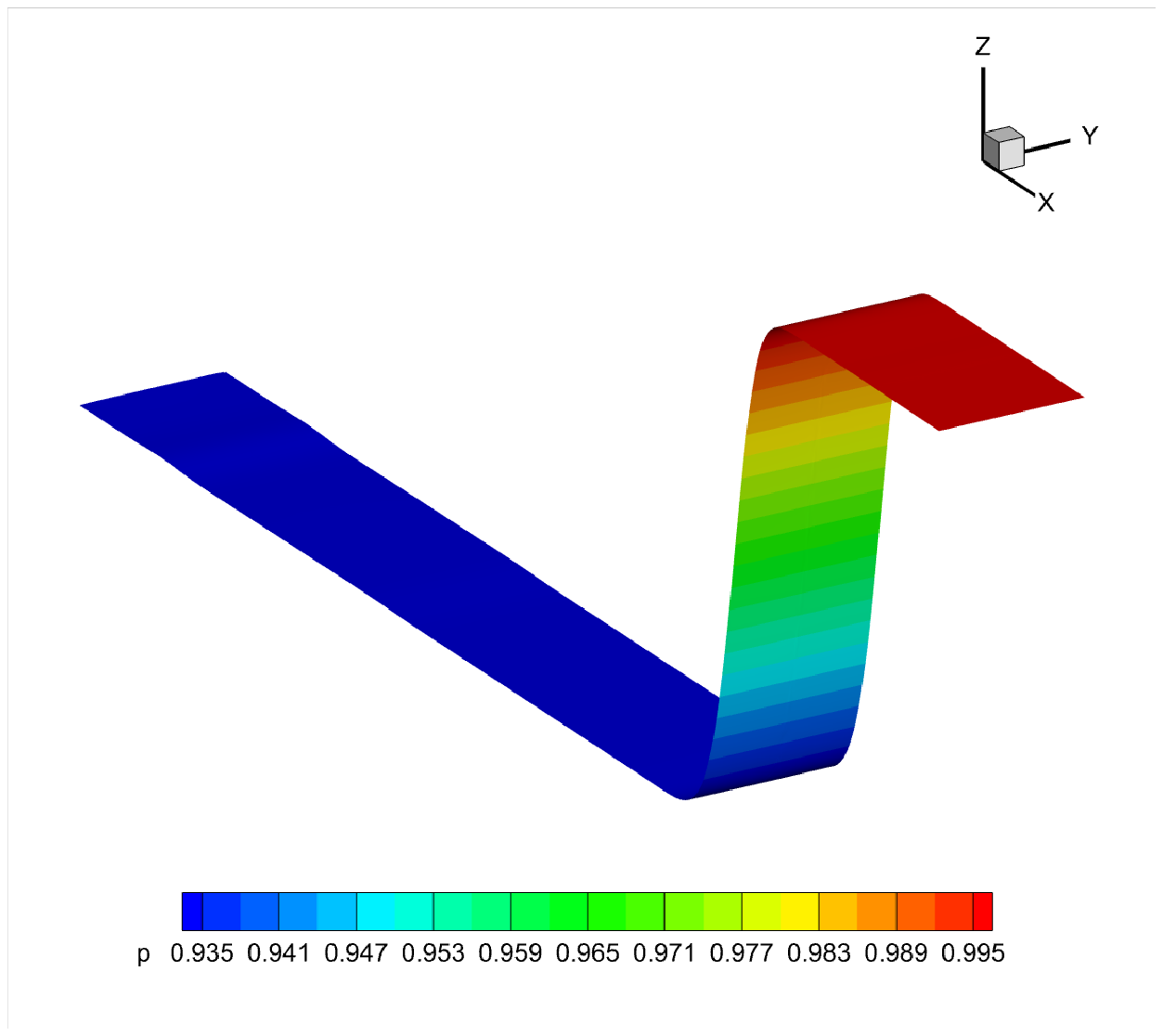} \\
		\end{tabular}
		\caption{Riemann problem RP1 at time $t=0.5$. Numerical results for density (top left), velocity (top right) and thermal impulse component $j_1$ (bottom left) compared against the reference solution. 1D cut at $y=0.05$ with 200 equidistant points is used for the SPFV solution. 3D view and color map of the pressure distribution are shown in the bottom right panel.}
		\label{fig.RP1}
	\end{center}
\end{figure}

\begin{figure}[!htbp]
	\begin{center}
		\begin{tabular}{cc}
			\includegraphics[trim= 5 5 5 5,clip,width=0.45\textwidth,keepaspectratio=true]{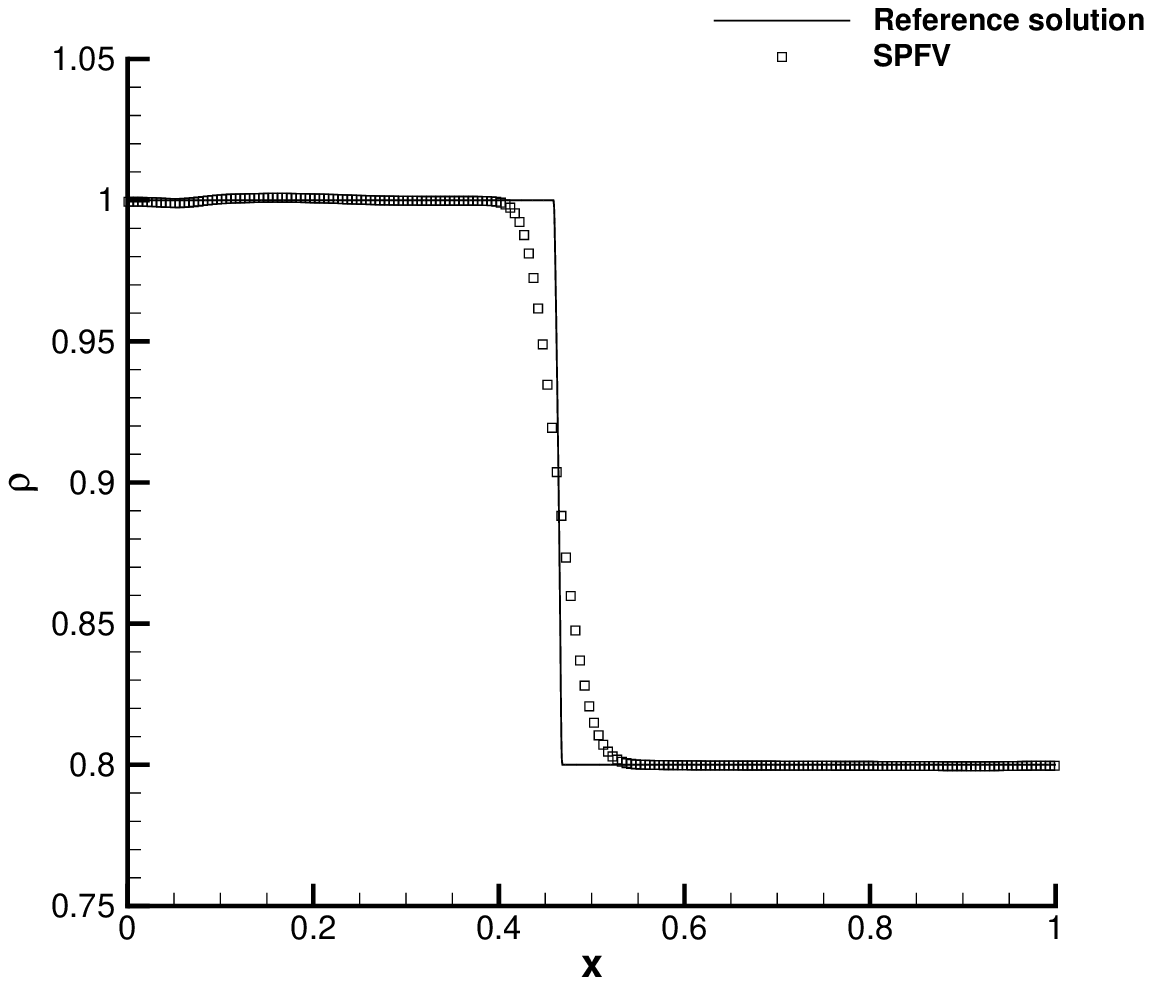} &
			\includegraphics[trim= 5 5 5 5,clip,width=0.45\textwidth,keepaspectratio=true]{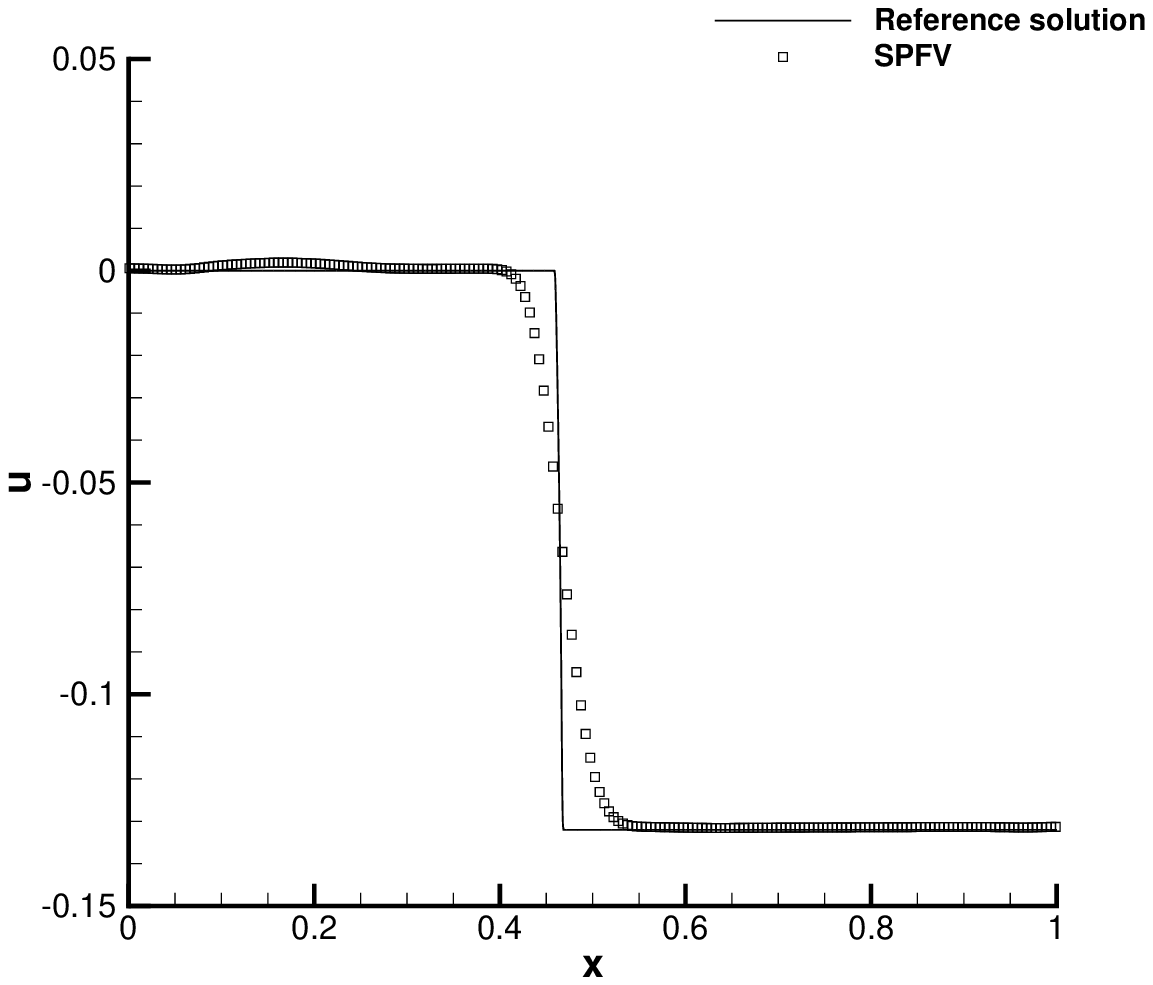} \\
			\includegraphics[trim= 5 5 5 5,clip,width=0.45\textwidth,keepaspectratio=true]{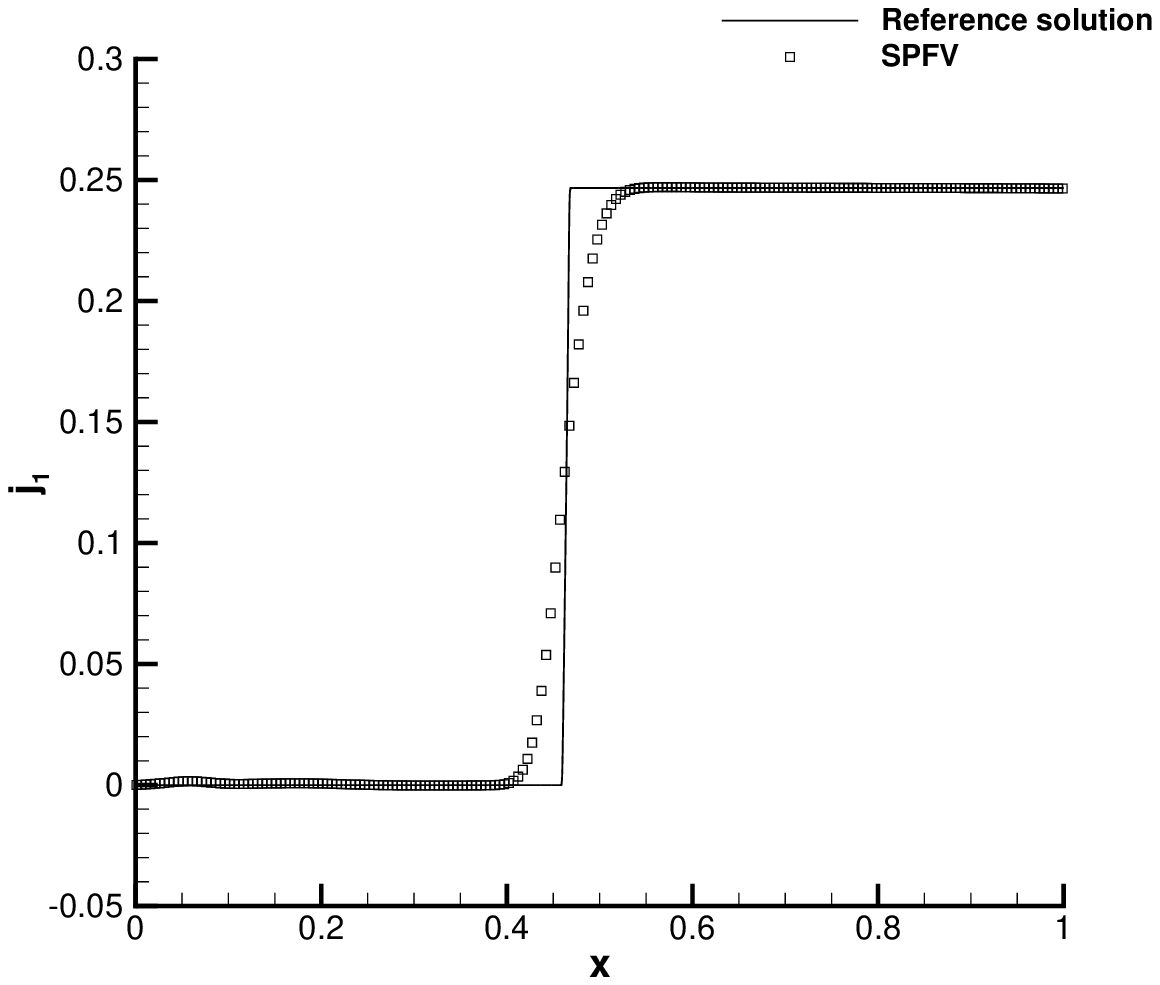} &
			\includegraphics[trim= 5 5 5 5,clip,width=0.45\textwidth,keepaspectratio=true]{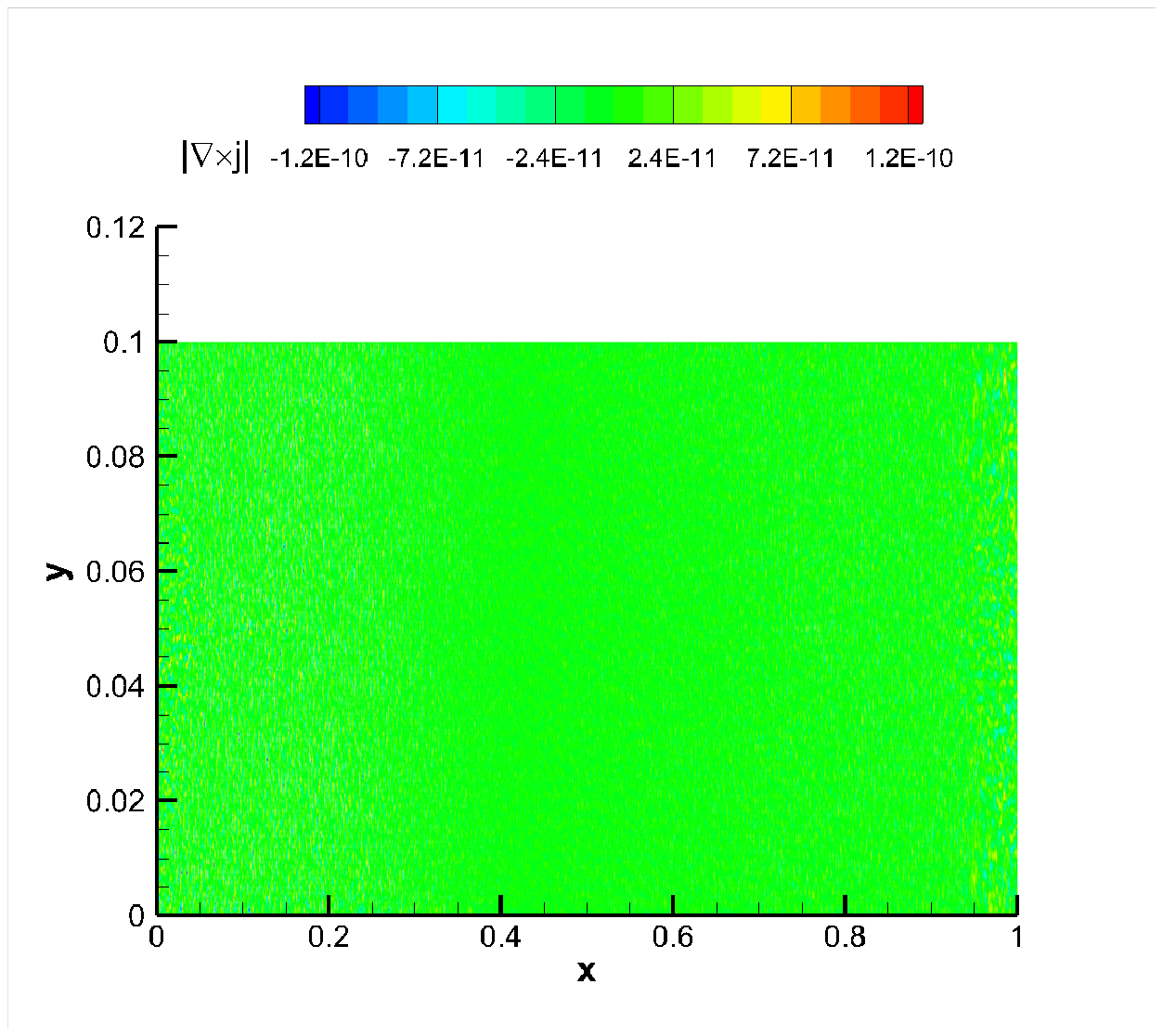} \\
		\end{tabular}
		\caption{Riemann problem RP2 at time $t=0.5$. Numerical results for density (top left), velocity (top right) and thermal impulse component $j_1$ (bottom left) compared against the reference solution. 1D cut at $y=0.05$ with 200 equidistant points is used for the SPFV solution. Color map of the absolute value of $\nabla \times \j$ is shown in the bottom right panel.}
		\label{fig.RP2}
	\end{center}
\end{figure}

To check that the new SPFV scheme can correctly capture the wave structures of the model, we run two Riemann problems (RP) proposed in \cite{dhaouadi2024eulerian}. The computational domain is $\Omega=[0,1]\times [0,0.1]$, and it is paved with a Voronoi tessellation of mesh size $h=1/1000$. We consider the homogeneous system, and the initial left and right state for RP1 read:
\begin{equation}
	\begin{array}{lcl}
		\mathbf{v}_L &=& \displaystyle \left( 0.8, \, -\frac{1}{4} \sqrt{1-\frac{\sqrt{13}}{8}}, \, 0, \, \frac{3}{4} + \frac{\sqrt{13}}{20}, \, \frac{1}{4} \sqrt{\frac{1}{15}(11+\sqrt{13})}, \, 0 \right), \\
		\mathbf{v}_R &=& \displaystyle \left( 1, \, 0, \, 0, \, 1, \, 0, \, 0 \right),
	\end{array}
\end{equation}
which are initially separated at $x_d=0.5$. For RP2, we simply switch the left and the right state, and the initial discontinuity is placed at $x_d=0.2$. The thermal impulse can be well initialized by choosing the piecewise linear thermal potential $\phi=j_1^{L,R} \, (x-x_d)$. RP1 considers an expansion shock, where the velocity is higher on the right side of the shock, while RP2 involves a compression fan solution. The final time of the simulations is $t_f=0.5$ and we set $\kappa=0.8$. Fig.~\ref{fig.RP1}-\ref{fig.RP2} show the numerical results which look in very good agreement with the reference, that is computed with a MUSCL-TVD scheme on a 1D mesh of $20'000$ cells.  Notice that, despite the unstructured character of the mesh size, the one-dimensional symmetry of the solution is very well preserved (see Fig.~\ref{fig.RP2}), as well as the curl-free property of $\j$ (see Fig.~\ref{fig.RP1}). Entropy is clearly dissipated for these test cases, thus $\as=0$ in \eqref{eqn.alpha}.

\subsection{Vortex problem with dissipation} \label{ssec.vortex}
\begin{figure}[!htbp]
	\begin{center}
		\begin{tabular}{cc}
			\includegraphics[trim= 5 5 5 5,clip,width=0.45\textwidth,keepaspectratio=true]{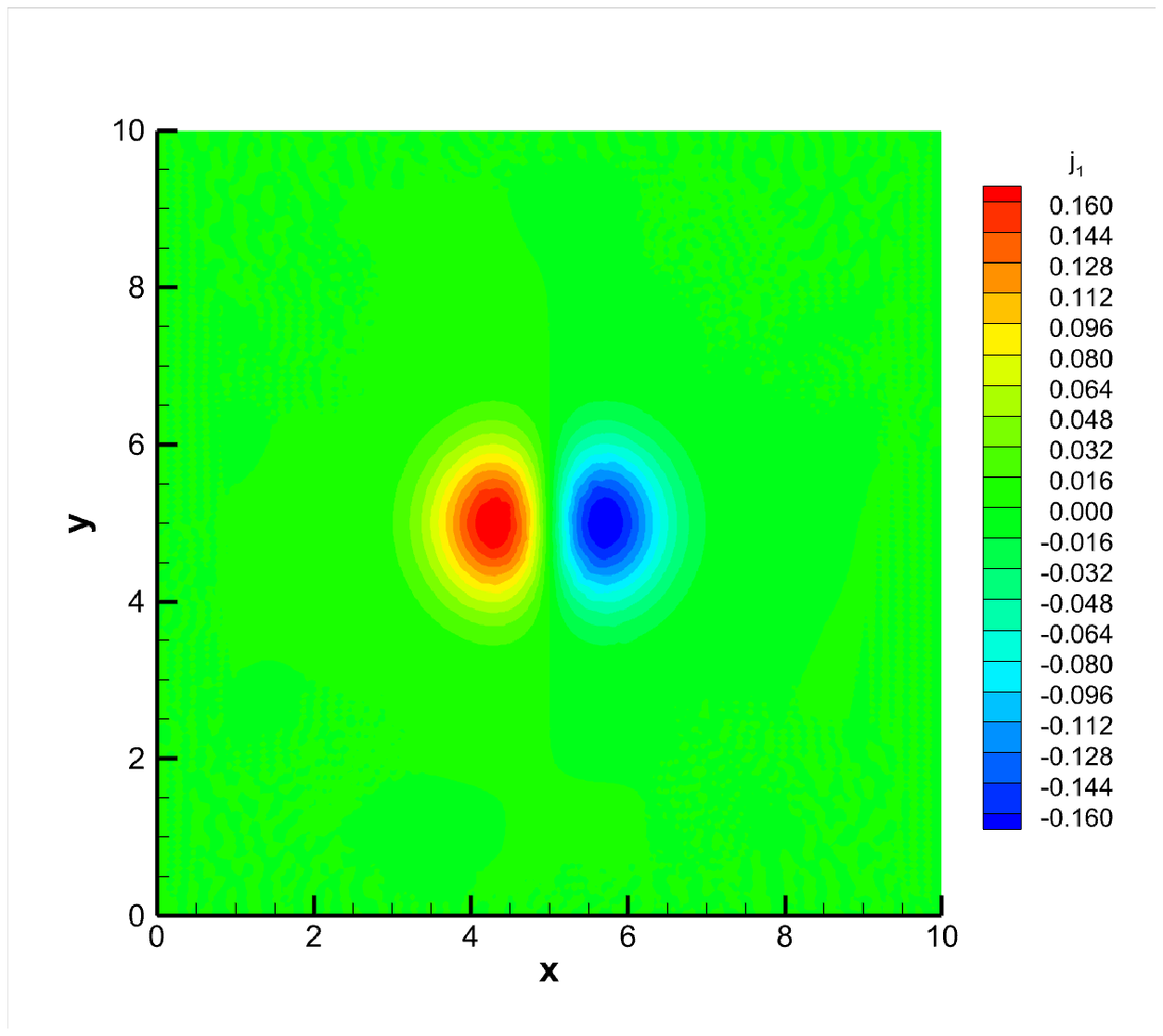} &
			\includegraphics[trim= 5 5 5 5,clip,width=0.45\textwidth,keepaspectratio=true]{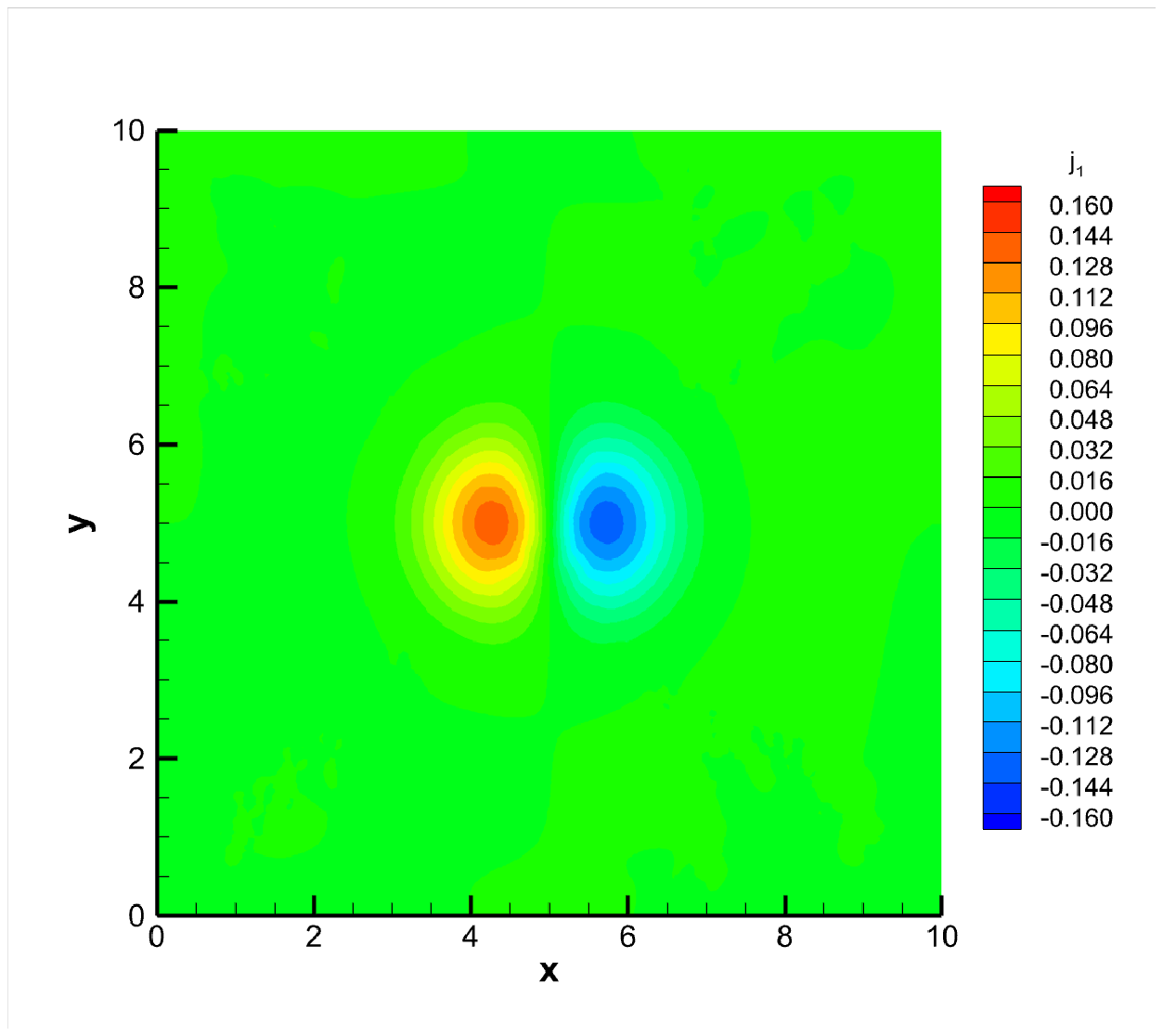} \\
			\includegraphics[trim= 5 5 5 5,clip,width=0.45\textwidth,keepaspectratio=true]{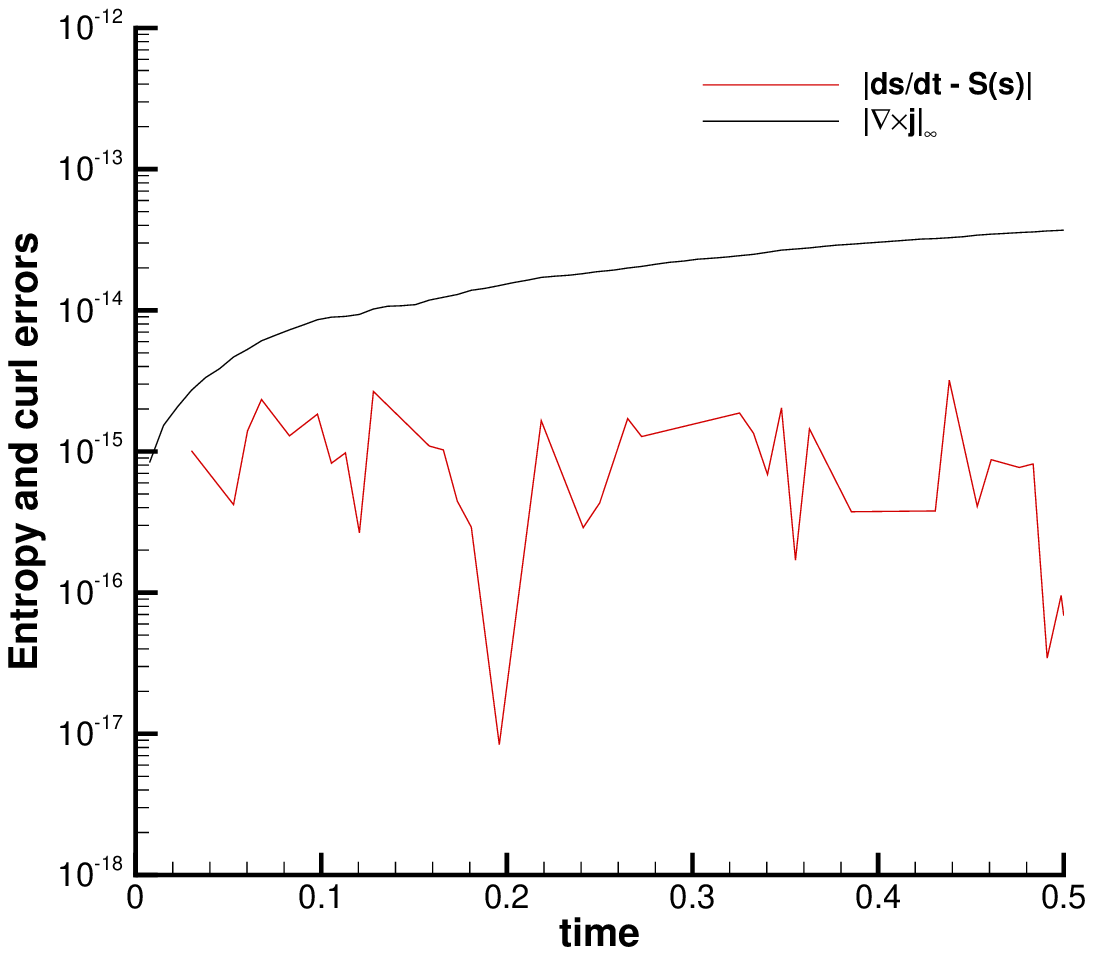} &
			\includegraphics[trim= 5 5 5 5,clip,width=0.45\textwidth,keepaspectratio=true]{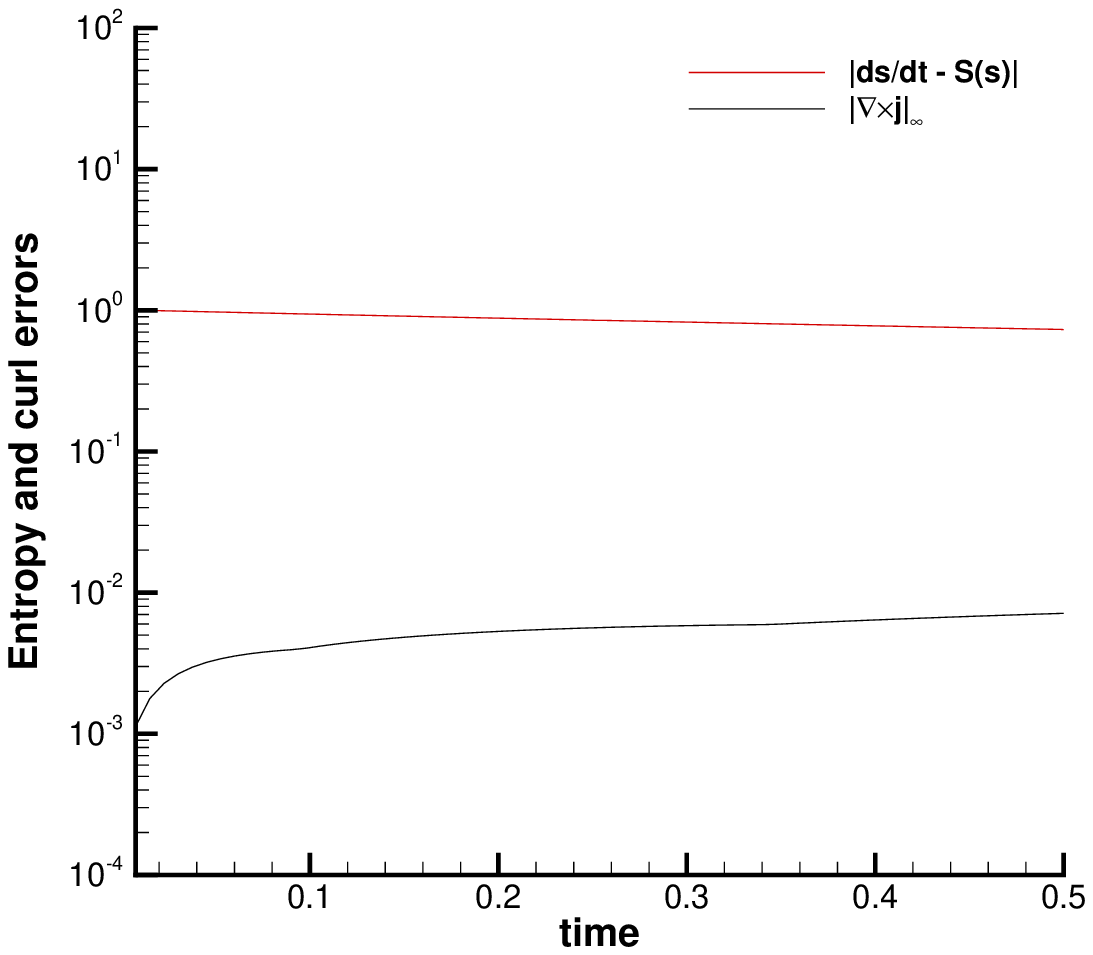} \\
		\end{tabular}
		\caption{Vortex problem at time $t=0.5$. Top: 21 equidistant contour lines for the thermal impulse component $j_1$. Bottom: time evolution of the error related to the curl of $\j$ and to the semi-discrete entropy equation with $S(s)=\kappa^2 \|\j\|^2 / (\rho \temp \tau)$. Left: curl preserving and thermodynamically compatible scheme (SPFV). Right: non compatible scheme.}
		\label{fig.vortex}
	\end{center}
\end{figure}

Here we aim at empirically showing the thermodynamic compatibility and the curl-preserving property. The computational domain $\Omega=[0,10]^2$ is discretized with Voronoi cells of size $h=1/10$. The isentropic vortex setup of \cite{HuShuTri} is used to initialize the Euler sub-system, that is
\begin{equation}
	\label{eqn.ICvortex}
	\rho = (1+\delta \temp)^{\frac{1}{\gamma-1}}, \quad \left( \begin{array}{c}
		u_1 \\ u_2 
	\end{array}\right) = \frac{5}{2\pi} e^{\frac{1-r^2}{2}} \left( \begin{array}{cc}
	5-y \\ x-5
\end{array}  \right), \quad \delta  \temp= -\frac{25 \, (\gamma-1)}{8 \gamma \pi^2}e^{1-r^2},
\end{equation}
with $p=\rho^{\gamma}$ and $r=\sqrt{x^2+y^2}$. Initially, $\phi=0$ thus $\j=0$. We set $\kappa=10^{-2}$ and $K=10^{-3}$, and the simulation is run until the final time $t_f=0.5$. The flow is no longer isentropic since dissipation is generated by the source term in the thermal impulse equation, thus we have to ensure that the semi-discrete entropy balance \eqref{eqn.sdeta} accounting for the source term $S(\rhos)=\kappa^2 \|\j\|^2/(\rho \temp \tau)$ is preserved. The relaxation time is kept constant and it is evaluated as $\tau=K/\kappa^2$, hence implying a curl-free evolution of the thermal impulse. The results are depicted in Fig.~\ref{fig.vortex} and compared with a non compatible scheme. The new curl preserving and thermodynamically compatible scheme is less dissipative while being stable, and it preserves very well the semi-discrete entropy equation and the curl-free behavior of $\j$.

\subsection{Explosion problem} \label{ssec.EP}
\begin{figure}[!htbp]
	\begin{center}
		\begin{tabular}{cc}
			\includegraphics[trim= 5 5 5 5,clip,width=0.45\textwidth,keepaspectratio=true]{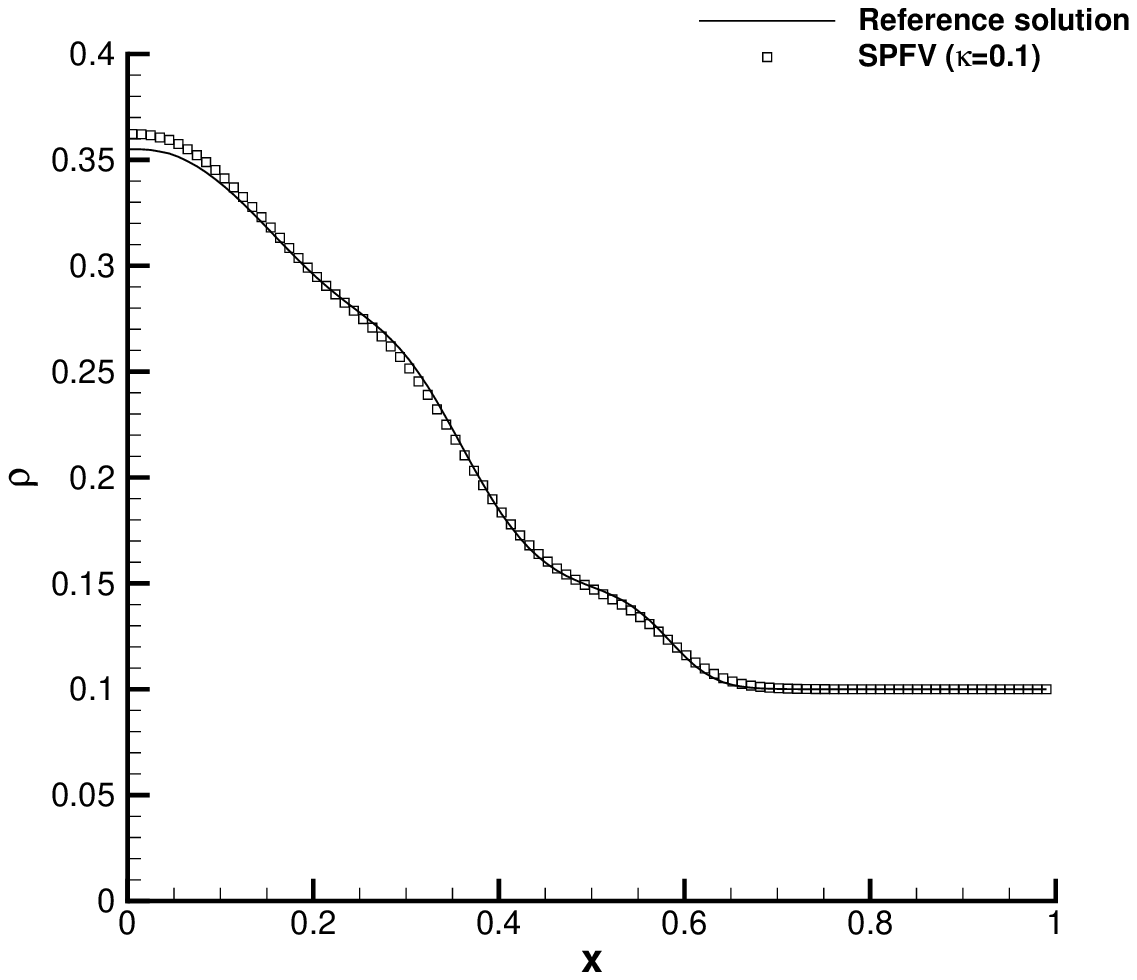} &
			\includegraphics[trim= 5 5 5 5,clip,width=0.45\textwidth,keepaspectratio=true]{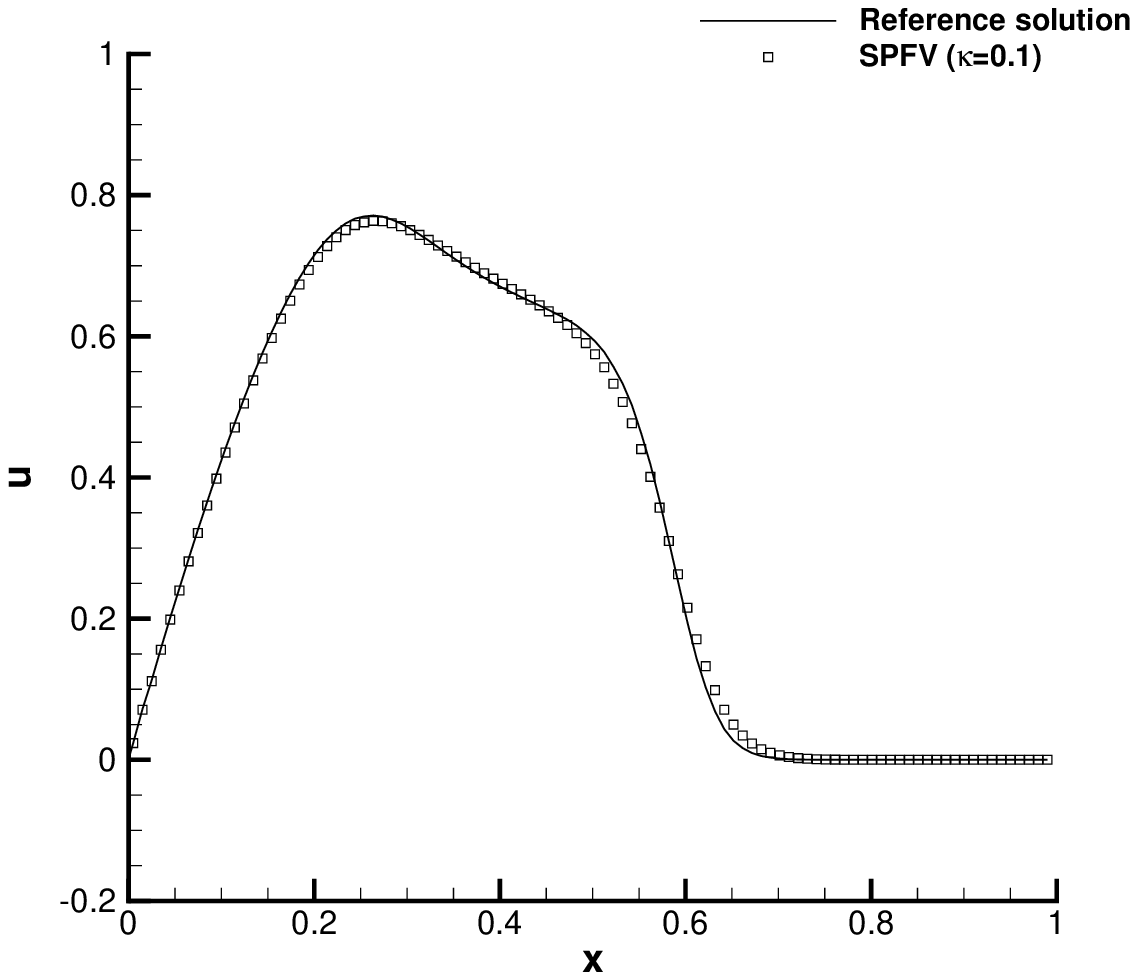} \\
			\includegraphics[trim= 5 5 5 5,clip,width=0.45\textwidth,keepaspectratio=true]{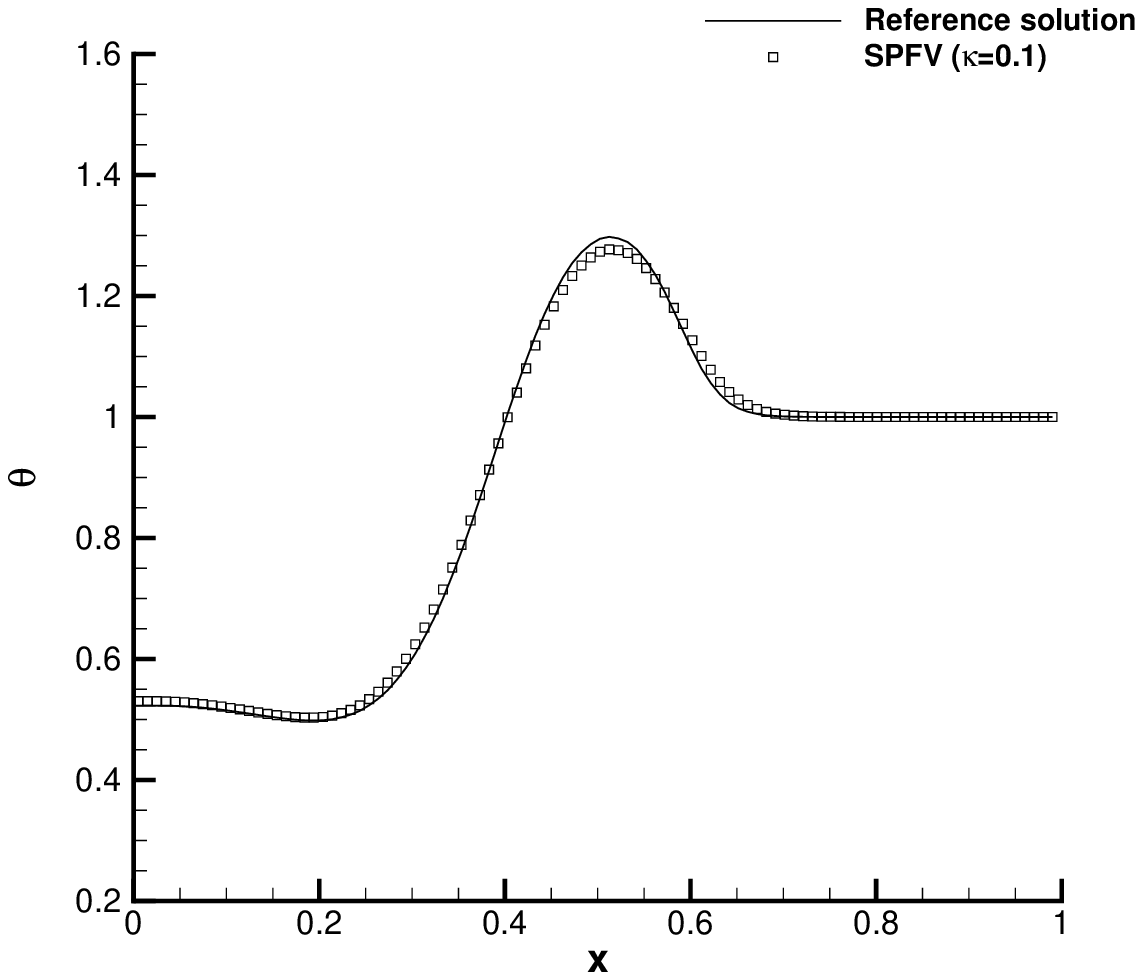} &
			\includegraphics[trim= 5 5 5 5,clip,width=0.45\textwidth,keepaspectratio=true]{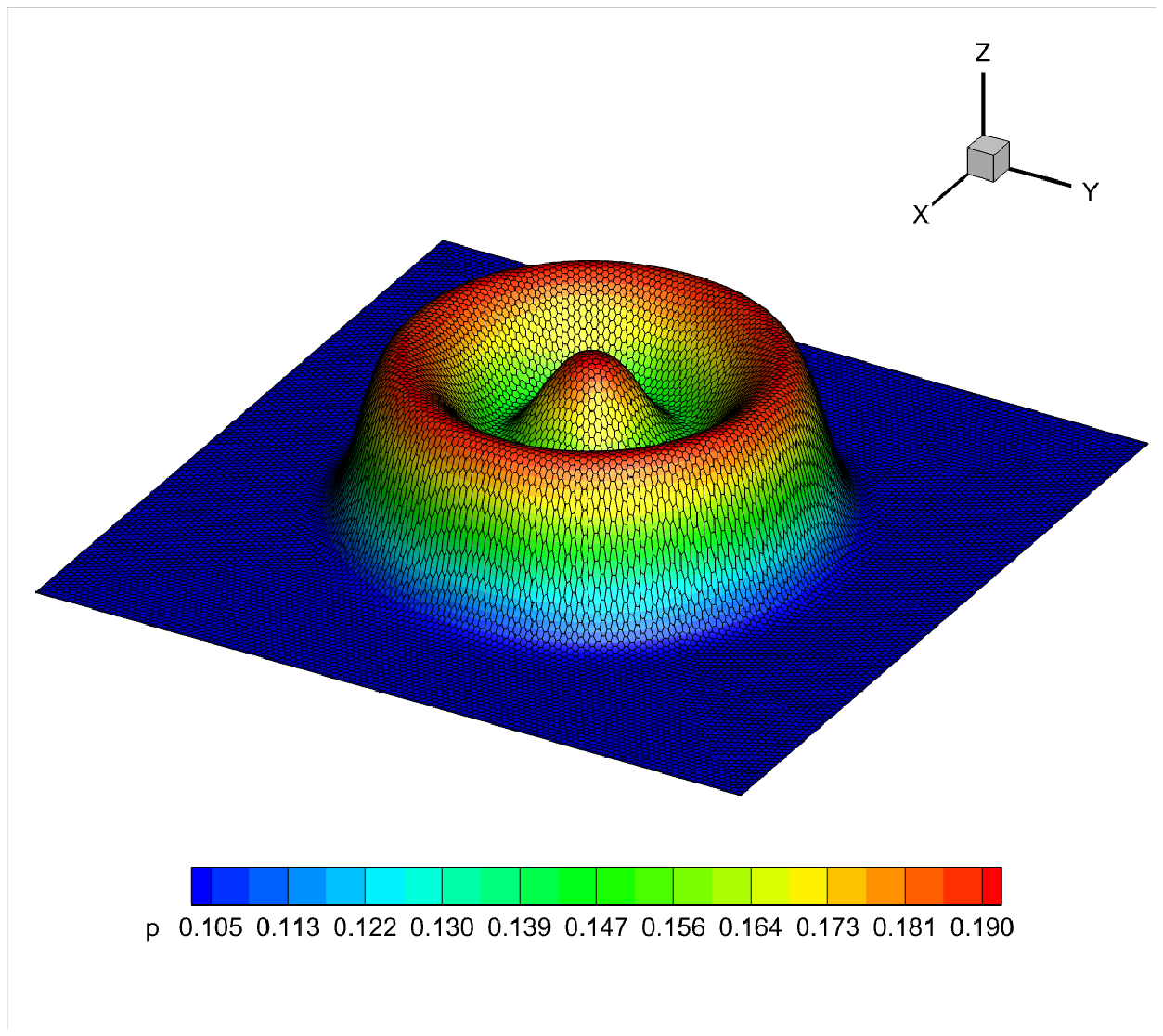} \\
		\end{tabular}
		\caption{Explosion problem at time $t=0.2$. Numerical results for density (top left), horizontal velocity (top right) and temperature (bottom left) compared against the reference solution. 1D cut at $y=0$ with 200 equidistant points is used for the SPFV solution. 3D view, Voronoi mesh and color map of the pressure distribution are shown in the bottom right panel.}
		\label{fig.EP2D}
	\end{center}
\end{figure}

Lastly, we show numerical evidence of the asymptotic preserving property of the SPFV scheme with the Euler equations supplemented with heat conduction. Let $\Omega=[-1,1]^2$ be the computational domain discretized by a Voronoi grid of size $h=1/60$. The initial condition is given by
\begin{equation}
	\begin{array}{lcll}
		\mathbf{v}_L &=& \displaystyle \left( 1, \, 0, \, 0, \, 1, \, 0, \, 0 \right), &\qquad \textnormal{for} \quad r<0.2, \\
		\mathbf{v}_R &=& \displaystyle \left( 0.1, \, 0, \, 0, \, 0.1, \, 0, \, 0 \right), &\qquad \textnormal{for} \quad r\geq 0.2,
	\end{array}
\end{equation}
with $r=\sqrt{x^2+y^2}$. Here, we take $\gamma=5/3$, $c_v=3/2$, and the thermal conductivity is $K=10^{-3}$. The dissipation source term is switched on, hence the relaxation parameter is $\kappa=0.1$, and the relaxation time is computed according to the asymptotic analysis relation \eqref{eqn.tau}, thus $\tau=(K \, \rho)/(\kappa^2 \, \theta)$. The reference solution is evaluated by solving the Euler system supplemented with Fourier heat conduction terms on the same computational grid using a MUSCL-TVD scheme until the time $t_f=0.2$. The results plot in Fig.~\ref{fig.EP2D} exhibits an overall excellent agreement of the SPFV scheme with the reference solution, showing the asymptotic preserving property of the new method. The cylindrical symmetry of the solution is also well preserved. Neither entropy nor curl-free preservation is obviously involved within this test case.

\section{Conclusions}
\label{sec:Conclusions}
In this paper we have presented a new finite volume scheme for the solution of the hyperbolic model of heat conducting flows \cite{dhaouadi2023structure}. The scheme makes use of a primal Voronoi grid and a dual triangular mesh to construct compatible curl-grad operators that ensure the curl-free property of any generic vertex-staggered vector field. We have derived a novel cell solver for the thermal impulse vector that is proven to yield conservation and local non-negative entropy production. An implicit discretization of the source term in the equation of the thermal impulse allows to retrieve, in the stiff relaxation limit, a consistent discretization of the Fourier law, hence making the scheme Asymptotic Preserving. Finally, we have proposed a novel subflux reconstruction technique to transfer the residuals from the dual to the primal mesh, so that thermodynamic compatibility can be ensured at the discrete level on the Voronoi grid. Eventually, thermodynamics compatibility can be achieved on general vertex-staggered grids. These properties have been proven, and the numerical results obtained for some basic test cases underline that the properties are also achieved in practice. We plan to extend these ideas to the GPR model of continuum mechanics \cite{PeshRom2014}, starting from the discretization recently forwarded in \cite{BoscheriGPRGCL}.

\bibliographystyle{siamplain}
\bibliography{references}

\end{document}